\pdfoutput=1
\documentclass[a4paper,11pt]{article}
\usepackage[english]{babel}
\usepackage[sc,osf]{mathpazo}
\usepackage[T1]{fontenc}
\usepackage[utf8]{inputenc}
\usepackage{amssymb}
\usepackage{amsthm}
\usepackage{thmtools}
\usepackage{graphicx}
\usepackage{amsmath}
\usepackage{enumerate}
\usepackage{pifont}
\usepackage{mathtools}
\usepackage{tensor}
\usepackage{stmaryrd}
\usepackage{placeins}
\usepackage{verbatim}
\usepackage{tikz-cd}
\usepackage{tikz}
\usetikzlibrary{shapes.geometric}
\usepackage{caption}
\usepackage{sidecap}
\usepackage{subcaption}
\usepackage{accents}
\usepackage{color}
\usepackage[pdfborder={0 0 0}]{hyperref}
\usepackage[nameinlink]{cleveref}  
\usepackage{imakeidx}
\usepackage{csquotes}
\usepackage{subfiles}
\usepackage[backend=biber,style=alphabetic,sorting=nty,giveninits]{biblatex}
\renewbibmacro{in:}{%
  \ifentrytype{article}{}{\printtext{\bibstring{in}\intitlepunct}}}
\tikzset{%
    symbol/.style={%
        draw=none,
        every to/.append style={%
            edge node={node [sloped, allow upside down, auto=false]{$#1$}}}
    }
}

\tikzset{
    rot90/.style={anchor=south, rotate=90, inner sep=.5mm}
}
\tikzset{
    rot320/.style={anchor=south, rotate=320, inner sep=.5mm}
}

\addto\extrasenglish{ 
 
 }

\declaretheorem[name=Theorem, numberwithin=section]{theorem}
\declaretheorem[name=Lemma, sibling=theorem]{lemma}
\declaretheorem[name=Proposition, sibling=theorem]{proposition}
\declaretheorem[name=Corollary, sibling=theorem]{corollary}

\declaretheorem[style=definition, name=Definition, sibling=theorem]{definition}
\declaretheorem[style=definition, name=Remark, sibling=theorem]{remark}
\declaretheorem[style=definition, name=Example, sibling=theorem]{example}

\declaretheorem[style=definition, name=Warning, sibling=theorem]{warrning}

\DeclareMathOperator{\Fun}{Fun}
\DeclareMathOperator*{\colim}{colim}
\DeclareMathOperator*{\holim}{holim}

\DeclareMathOperator{\id}{id}

\DeclareMathOperator{\Pro}{Pro}
\DeclareMathOperator{\Ind}{Ind}
\DeclareMathOperator{\ev}{ev}

\DeclareMathOperator{\Hom}{Hom}
\DeclareMathOperator{\Map}{Map}
\DeclareMathOperator{\map}{map}
\DeclareMathOperator{\Ar}{Ar}
\DeclareMathOperator{\cof}{cof}
\DeclareMathOperator{\we}{we}
\DeclareMathOperator{\fib}{fib}
\DeclareMathOperator{\cosk}{cosk}
\DeclareMathOperator{\sk}{sk}
\DeclareMathOperator{\Ho}{Ho}
\DeclareMathOperator{\Sp}{Sp}
\DeclareMathOperator{\Sing}{Sing_\mathit{J}}
\DeclareMathOperator{\SingTop}{Sing}
\newcommand{\wh}{\widehat}
\newcommand{\wt}{\widetilde}

\newcommand{\ul}{\underline}
\newcommand{\ultimes}{\mathbin{\ooalign{$\hidewidth\underline{\times}\hidewidth$\cr$\phantom{\times}$}}}
\newcommand{\bbN}{\mathbb N}

\newcommand{\bfC}{\mathbf C}

\newcommand{\bfE}{\mathcal E}
\newcommand{\bfT}{\mathbf T}
\newcommand{\bfG}{\mathbf G}
\newcommand{\bftwo}{\mathbf 2}
\newcommand{\s}{\mathbf s}
\newcommand{\Fin}{\mathbf{Fin}}
\newcommand{\Set}{\mathbf{Set}}
\newcommand{\Cat}{\mathbf{Cat}}
\newcommand{\Grp}{\mathbf{Grp}}
\newcommand{\Grpd}{\mathbf{Grpd}}
\newcommand{\Stone}{\mathbf{Stone}}
\newcommand{\Top}{\mathbf{Top}}
\newcommand{\FinsSet}{\mathbf{sSet}_\mathrm{fin}}
\newcommand{\LsSet}{\mathbf{L}}

\newcommand{\bisSet}{\mathbf{bisSet}}
\newcommand{\biswhSet}{\mathbf{bis}\wh\Set}

\makeatletter
\newcommand{\oset}[3][0ex]{%
  \mathrel{\mathop{#3}\limits^{
    \vbox to#1{\kern-2\ex@
    \hbox{$\scriptstyle#2$}\vss}}}}
\makeatother

\newcommand{\cofarrow}{\rightarrowtail}
\newcommand{\trivcofarrow}{\oset[-.5ex]{\sim}{\rightarrowtail}}
\newcommand{\wearrow}{\oset[-.16ex]{\sim}{\longrightarrow}}
\newcommand{\fibarrow}{\twoheadrightarrow}
\newcommand{\trivfibarrow}{\oset[-.5ex]{\sim}{\twoheadrightarrow}}

\renewcommand{\subset}{\subseteq}

\mathchardef\mhyphen="2D

\title{Simplicial model structures on pro-categories}
\author{Thomas Blom and Ieke Moerdijk}
\date{\today}

\begin{document}

\maketitle
\begin{abstract}
    We describe a method for constructing simplicial model structures on ind- and pro-categories. Our method is particularly useful for constructing ``profinite'' analogues of known model categories. Our construction quickly recovers Morel's model structure for pro-$p$ spaces and Quick's model structure for profinite spaces, but we will show that it can also be applied to construct many interesting new model structures. In addition, we study some general properties of our method, such as its functorial behaviour and its relation to Bousfield localization. We compare our construction to the $\infty$-categorical approach to ind- and pro-categories in an appendix.
\end{abstract}

\tableofcontents


\section{Introduction}

In \cite{Quick2008Profinite,Quick2011Continuous}, Quick constructed a fibrantly generated Quillen model structure on the category of simplicial profinite sets that models the homotopy theory of ``profinite spaces''.
This can be seen as a continuation of Morel's work in \cite{Morel1996EnsemblesProfinis}, where, for a given prime $p$, he presented a model structure on the same category that models the homotopy theory of ``pro-$p$ spaces''.

The purpose of this paper is to present a new and uniform method that immediately gives these two model structures, as well as many others. For example, while Quick's
model structure is in a sense derived from the classical homotopy theory
of simplicial sets, our method also applies to the Joyal model
structure, thus providing a homotopy theory of profinite $\infty$-categories. Our construction can also be used to obtain a model category of profinite $P$-stratified spaces, where $P$ is a finite poset, whose underlying $\infty$-category is the $\infty$-category of profinite $P$-stratified spaces defined in \cite{BarwickGlasmanHaine2018ExodromyV7}.

One general form that our results take is the following version of
pro-completion of model categories:

\begin{theorem}\label{thm:ApproximateTheorem}
Let $\bfE$ be a simplicial model category in which every object is cofibrant and let $\bfC$ be an (essentially) small full subcategory of $\bfE$ closed under finite limits and cotensors by finite simplicial sets. Then for any collection $\bfT$ of fibrant objects in $\bfC$, the pro-completion $\Pro(\bfC)$ carries a fibrantly generated simplicial model structure with the following properties:
\begin{enumerate}[(1)]
    \item \label{thm:ApproximateTheoremItem1} The weak equivalences are the $\bfT$-local equivalences; that is, $f \colon C \to D$ is a weak equivalence if and only if 
    \[f^* \colon \Map(D,t) \to \Map(C,t) \]
    is a weak equivalence of simplicial sets for any $t \in \bfT$.
    \item \label{thm:ApproximateTheoremItem2} Every object in $\Pro(\bfC)$ is again cofibrant.
    \item \label{thm:ApproximateTheoremItem3} The inclusion $\bfC \hookrightarrow \bfE$ induces a simplicial Quillen adjunction $\bfE \rightleftarrows \Pro(\bfC)$.
    \item \label{thm:ApproximateTheoremItem4}If $\bfT \subset \bfC$ is closed under pullbacks along fibrations and cotensors by finite simplicial sets, then the underlying $\infty$-category of this model structure on $\Pro(\bfC)$ is equivalent to $\Pro(N(\bfT))$, where $N(\bfT)$ denotes the homotopy coherent nerve of the full simplicial subcategory of $\bfE$ spanned by the objects of $\bfT$.
\end{enumerate}
\end{theorem}

The model structures of Quick and Morel mentioned above can be obtained from this theorem by appropriately choosing a full subcategory $\bfC$ of $\s\Set$ and a collection $\bfT$ of fibrant objects. Another known model structure that can be recovered from the above theorem is the model structure for ``profinite groupoids'' constructed by Horel in \cite[\S 4]{Horel2017ProfiniteOperads}.

The new model category $\Pro(\bfC)$ is a kind of pro-completion of $\bfE$ with respect to the pair $(\bfC,\bfT)$, and could be denoted $\wh \bfE$ or $\bfE^\wedge_{(\bfC,\bfT)}$. The left adjoint $\bfE \to \Pro(\bfC)$ of the Quillen adjunction mentioned in item \eqref{thm:ApproximateTheoremItem3} can be seen as a ``pro-$\bfC$ completion'' functor. For the model structures of Morel, Quick and Horel mentioned above, this functor agrees with the profinite completion functor.

We would like to point out that the above formulation is slightly incomplete since there are multiple ways of choosing sets of generating (trivial) fibrations, which theoretically could lead to different model structures on $\Pro(\bfC)$, though always with the weak equivalences as described above. A noteworthy fact is that the above theorem also holds for model categories enriched over the Joyal model structure on simplicial sets, so in particular it applies to the Joyal model structure itself. In this case, the model structure obtained on $\Pro(\bfC)$ is enriched over the Joyal model structure, but not necessarily over the classical Kan-Quillen model structure on $\s\Set$. Another fact worth mentioning is that there exist many simplicial model categories satisfying the hypotheses of the above theorem; that is, all objects being cofibrant. Indeed, by a result of Dugger (Corollary 1.2 of \cite{Dugger2001CombinatorialPresentations}), any combinatorial model category is Quillen equivalent to such a simplicial model category.

Even though we are mostly interested in model structures on pro-categories, we will first describe our construction in the context of ind-categories, and then dualize those results. We have chosen this approach since in the case of ind-categories our construction produces cofibrantly generated model categories, which to most readers will be more familiar territory than that of fibrantly generated model categories. In addition, this will make it clear that the core of our argument, which is contained in \autoref{sec:IndCompletedModelStructure}, only takes a few pages. Another reason for describing our construction in the context of ind-categories is that an interesting example occurs there: if we apply our construction to a well-chosen full subcategory of the category of topological spaces, then we obtain a model category that is Quillen equivalent to the usual Quillen model structure on $\Top$, but that has many favourable properties, such as being combinatorial.

Our original motivation partly came from the desire to have a full fledged Quillen style homotopy theory of profinite $\infty$-operads, by using the category of dendroidal Stone spaces (i.e.~dendroidal profinite sets). However, not every object is cofibrant in the operadic model structure for dendroidal sets, so the methods from the current paper do not apply directly to this case. The extra work needed to deal with objects that are not cofibrant is of a technical nature, and very specific to the example of dendroidal sets. For this reason, we have decided to present this case separately, see \cite{BlomMoerdijk2021ProfiniteOperadsArxivV1}.

\paragraph{Relation to the construction by Barnea-Schlank.} There are several results in the literature that describe general methods for constructing model structures on ind- or pro-categories. The construction in the current paper is quite close in spirit to that by Barnea and Schlank in \cite{BarneaSchlank2016ProSimplicialSheaves}. They show that if $\bfC$ is a category endowed with the structure of a ``weak fibration category'', then there exists an ``induced'' model structure on $\Pro(\bfC)$ provided some additional technical requirements are satisfied. However, there are important examples of model structures on pro-categories that are not of this form. For example, Quick's model structure is not of this kind, as explained just above Proposition 7.4.2 in \cite{BarneaHarpazHorel2017}. In the present paper, we prove the existence of a certain model structure on the pro-category of a simplicial category endowed with the extra structure of a so-called ``fibration test category'' (defined in \autoref{definition:FibTestCat}). While the definition of a fibration test category given here seems less general than that of a weak fibration category, there are many interesting examples where it is easy to prove that a category is a fibration test category while it is not clear whether this category is a weak fibration category in the sense of \cite{BarneaSchlank2016ProSimplicialSheaves}. In particular, Quick's model structure can be obtained through our construction, see \autoref{example:KanQuillenLeanFibrationTestCategory} and \autoref{corollary:QuicksModelStructureAgrees}. Another advantage is that we do not have to check the technical requirement of ``pro-admissibility'' (see Definition 4.4 of \cite{BarneaSchlank2016ProSimplicialSheaves}) to obtain a model structure on $\Pro(\bfC)$, which is generally not an easy task. We also believe that our description of the weak equivalences in $\Pro(\bfC)$, namely as the $\bfT$-local equivalences for some collection of objects $\bfT$, is often more natural and flexible than the one given in \cite{BarneaSchlank2016ProSimplicialSheaves}. It is worth pointing out that if both our model structure and that of \cite{BarneaSchlank2016ProSimplicialSheaves} on $\Pro(\bfC)$ exist, then they agree by \autoref{remark:BarneaSchlankVsFibrationTestCategory} below.

\paragraph{Overview of the paper.} In \autoref{sec:Preliminaries}, we will establish some terminology and mention a few facts on simplicial model categories and ind- and pro-categories. We will then describe our general construction of the model structure for ind-categories in \autoref{sec:IndCompletedModelStructure}.  We illustrate our construction with an example in \autoref{sec:ConvenientCategorySpaces}, where we construct a convenient model category of spaces. In \autoref{sec:ProCompletedModelStructure}, we dualize our results to the context of pro-categories, and illustrate this dual construction with many examples. We show that some of these examples coincide with model structures that are already known to exist in \autoref{sec:Comparison}, such as Quick's and Morel's model structures. We then continue the study of our construction in \autoref{sec:QuillenPairs}, where we discuss its functorial behaviour, and in \autoref{sec:BousfieldLocalizations}, where we prove results about the existence of certain Bousfield localizations. The latter section also contains the proof of \autoref{thm:ApproximateTheorem}, except for item \eqref{thm:ApproximateTheoremItem4}. We then give a detailed discussion of two examples in \autoref{sec:CompleteSegalProfiniteVSProfiniteQuasiCats}; namely the model structure for complete Segal profinite spaces and the model structure for profinite quasi-categories. In the \hyperref[appendix:InftyApproach]{appendix}, we compare our construction to the $\infty$-categorical approach to ind- and pro-categories.

\paragraph{Acknowledgements} We would like to thank the referees for numerous comments that helped improve the exposition.


\section{Preliminaries}\label{sec:Preliminaries}
In this section, we will briefly review some basic definitions concerning simplicial objects, and then discuss ind- and pro-categories.

\subsection{Simplicial conventions}\label{ssec:SimplicialConventions}

We assume the reader to be familiar with the basic theory of simplicial sets \cite{May1967SimplicialObjects,Lamotke1968Semisimpliziale,GabrielZisman1967Calculus,GoerssJardine2009Simplicial}. We will say that a simplicial set $X$ is \emph{skeletal} if it is $n$-skeletal for some natural number $n$, i.e.~if the map $\sk_n X \to X$ is an isomorphism. Dually, $X$ is \emph{coskeletal} if $X \to \cosk_n X$ is an isomorphism for some $n$. Recall that a simplicial set $X$ is \emph{degreewise finite} if each $X_n$ is a finite set, and \emph{finite} if it has finitely many non-degenerate simplices. Note that the latter is equivalent to $X$ being degreewise finite and skeletal. We will say that a simplicial set is \emph{lean} if it is degreewise finite and coskeletal, and write $\LsSet$ for the full subcategory of $\s\Set$ on the lean simplicial sets. One can show that if $X$ is a lean simplicial set and if $Y$ is a degreewise finite simplicial set, then the cotensor $X^Y = \Map(Y,X)$ is again a lean simplicial set.

Most categories we deal with are \emph{simplicial categories}, i.e.~categories enriched over simplicial sets. Moreover, they will generally be required to have tensors or cotensors by finite simplicial sets. For objects $c$ and $d$ in a simplicial category $\bfC$, we will write $\Map(c,d)$ for the simplicial hom set. Recall that for a morphism $c \to d$ in $\bfC$ and a morphism $U \to V$ of simplicial sets, the \emph{pushout-product map} is the map
\[d \otimes U \cup_{c \otimes U} c \otimes V \to d \otimes V ,\]
which makes sense in $\bfC$ if the necessary pushouts and tensors exist. Dually, we refer to 
\[c^V \to c^U \times_{d^U} d^V \]
as the \emph{pullback-power map} (if it exists). If given another morphism $a \to b$ in $\bfC$, we refer to
\[\Map(b,c) \to \Map(a,c) \times_{\Map(a,d)} \Map(b,d) \]
as a \emph{pullback-power map} as well. Note that this map always exist.

We assume the reader to be familiar with the basic theory of Quillen model categories \cite{Hovey1999ModelCats,Hirschhorn2003Model}. Basic examples include the classical Kan-Quillen model structure on simplicial sets which we denote by $\s\Set_{KQ}$, and the Joyal model structure $\s\Set_J$ modelling the homotopy theory of $\infty$-categories \cite{Lurie2009HTT}. A \emph{simplicial model category} is a model category $\bfE$ that is enriched, tensored and cotensored over simplicial sets, and that satisfies the additional axiom SM7 phrased in terms of pullback-power maps, or dually in terms of pushout-product maps (see e.g. Definition 9.1.6 and Proposition 9.3.7 of \cite{Hirschhorn2003Model} or Definition II.2.2 of \cite{Quillen1967HomotopicalAlgebra}). We emphasize that we will use this terminology in a somewhat non-standard way. Namely, by a simplicial model category, we will either mean that the axiom SM7 holds with respect to the Kan-Quillen model structure or the Joyal model structure. Whenever it is necessary to emphasize the distinction, we will call a simplicial model category of the former kind a \emph{$\s\Set_{KQ}$-enriched model category} and the latter a \emph{$\s\Set_J$-enriched model category}. Note that any $\s\Set_{KQ}$-enriched model category is $\s\Set_J$-enriched, since $\s\Set_{KQ}$ is a left Bousfield localization of $\s\Set_J$.

We will make use of the following fact about the (categorical) fibrations in $\s\Set_J$.

\begin{lemma}\label{lem:SkeletalJ}
There exists a set $M$ of maps between finite simplicial sets such that a map between quasi-categories $X \to Y$ is a fibration in $\s\Set_J$ if and only if it has the right lifting property with respect to all maps in $M$.
\end{lemma}

\begin{proof}
Let $H$ denote the simplicial set obtained by gluing two $2$-simplices to each other along the edges opposite to the $0$th and $2$nd vertex, respectively, and then collapsing the edges opposite to the $1$st vertex to a point in both of these $2$-simplices. This means that $H$ looks as follows, where the dashed lines represent the collapsed edges:
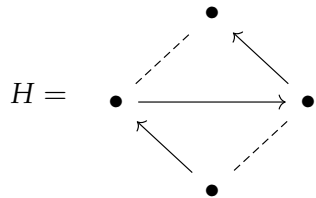
\[\begin{array}{rl}
   H =  &  \begin{tikzcd}[column sep = scriptsize]
 & \bullet & \\
\bullet \ar[ur, dashed, no head] \ar[rr] & & \bullet \ar[ul] \ar[dl, dashed, no head] \\
 & \bullet \ar[ul] &
\end{tikzcd}
\end{array} \]
A map from $H$ into a quasicategory $X$ consists of an arrow $f \in X_1$, a left and right homotopy inverse $g,h \in X_1$ and homotopies $gf \sim \id$ and $fh \sim \id$. Let $\{0\} \hookrightarrow H$ denote the inclusion of the leftmost vertex into $H$. It follows from Corollary 2.4.6.5 of \cite{Lurie2009HTT} that if $X \to Y$ is an inner fibration between quasicategories that has the right lifting property with respect to $\{0\} \hookrightarrow H$, then it is an categorical fibration. The converse is also true. To see this, note that for any quasicategory $Z$, a map $H \to Z$ lands in the largest Kan complex $k(Z)$ contained in $Z$. Since $\{0\} \hookrightarrow H$ is a weak homotopy equivalence, we see that $\Map(H,Z) = \Map(H,k(Z)) \simeq \Map(\{0\},k(Z)) = \Map(\{0\},Z)$, so the inclusion $\{0\} \hookrightarrow H$ is a categorical equivalence. In particular, any categorical fibration has the right lifting property with respect to $\{0\} \to H$. We conclude that the set $M = \{\Lambda^n_k \hookrightarrow \Delta^n \mid 0 < k < n \} \cup \{\{0\} \hookrightarrow H\}$ has the desired properties.
\end{proof}

\subsection{Ind- and pro-categories}\label{ssec:ProAndIndCategories}

In this section we recall some basic definitions concerning ind- and pro-categories. Most of these will be familiar to the reader, with the possible exception of \autoref{thm:IndValuedPresheavesIsIndSkeletalPresheaves} below. For details, we refer the reader to \cite[Exposé 1]{grothendieck1972theorie}, \cite[\S 2.1]{EdwardsHastings1976CechSteenrod}, \cite[Appendix]{ArtinMazur1986Etale} and \cite{Isaksen2002Calculating}. In the discussion below, all (co)limits are asssumed to be \textbf{small}.

For a category $\bfC$, its ind-completion $\Ind(\bfC)$ is obtained by freely adjoining filtered (or directed) colimits to $\bfC$. Dually, the free completion under cofiltered limits is denoted $\Pro(\bfC)$. This in particular means that $\Pro(\bfC)^{op} = \Ind(\bfC^{op})$, so any statement about ind-categories dualizes to a statement about pro-categories and vice versa. We will therefore mainly discuss ind-categories here and leave it to the reader to dualize the discussion.

One way to make the above precise, is to define the objects in $\Ind(\bfC)$ to be all diagrams $I \to \bfC$ for all filtered categories $I$. Such objects are called \emph{ind-objects} and denoted $C = \{c_i\}_{i \in I}$. The morphisms between two such objects $C = \{c_i\}_{i \in I}$ and $D = \{d_j\}_{j \in J}$ are defined by
\begin{equation}\label{equation:DefinitionHomSetIndCategory}
    \Hom_{\Ind(\bfC)}(C,D) = \lim_i \colim_j \Hom_{\bfC}(c_i,d_j).
\end{equation}
If $\bfC$ is a simplicial category, then $\Ind(\bfC)$ can be seen as a simplicial category as well. The enrichment is expressed by a formula similar to \eqref{equation:DefinitionHomSetIndCategory}, namely
\[\Map(\{c_i\}, \{d_j\}) = \lim_i \colim_j \Map(c_i, d_j).\]
One can define the pro-category $\Pro(\bfC)$ of a (simplicial) category $\bfC$ as the category of all diagrams $I \to \bfC$ for all cofiltered $I$, and with (simplicial) hom sets dual to the ones above. An object in $\Pro(\bfC)$ is called a \emph{pro-object}. One could also simply define $\Pro(\bfC)$ as $\Ind(\bfC^{op})^{op}$.

It can be shown that any object in $\Ind(\bfC)$ is isomorphic to one where the indexing category $I$ is a directed poset, and dually that any object in $\Pro(\bfC)$ is isomorphic to one that is indexed by a codirected poset (see Proposition 8.1.6 of \cite[Exposé 1]{grothendieck1972theorie}, or Theorem 2.1.6 of \cite{EdwardsHastings1976CechSteenrod} with a correction just after Corollary 3.11 of \cite{BarneaSchlank2015NewModel}).

There is a fully faithful embedding $\bfC \hookrightarrow \Ind(\bfC)$ sending an object $c$ to the constant diagram with value $c$, again denoted $c$. We will generally identify $\bfC$ with its image in $\Ind(\bfC)$ under this embedding. This embedding preserves all limits and all finite colimits that exist in $\bfC$. The universal property of $\Ind(\bfC)$ states that $\Ind(\bfC)$ has all filtered colimits and that any functor $F \colon \bfC \to \bfE$, where $\bfE$ is a category that has all filtered colimits, has an essentially unique extension to a functor $\wt F \colon \Ind(\bfC) \to \bfE$ that preserves filtered colimits. This extension can be defined explicitly by $\wt F(\{c_i\}) = \colim_i F(c_i)$.

Recall that if $\bfE$ is a category that has all filtered colimits, then an object $c$ in $\bfE$ is called \emph{compact} if $\Hom_\bfE(c,-)$ commutes with filtered colimits. The dual notion is called \emph{cocompact}. One can deduce from the definition of the morphisms in $\Ind(\bfC)$ that any object in the image of $\bfC \hookrightarrow \Ind(\bfC)$ is compact. Dually, the objects of $\bfC$ are cocompact in $\Pro(\bfC)$.

There is the following recognition principle for ind-completions, of which we leave the proof to the reader.

\begin{lemma}[Recognition principle]\label{lemma:RecognitionPrincipleInd}
Let $\bfE$ be a category closed under filtered colimits and let $\bfC \hookrightarrow \bfE$ be a full subcategory. If
\begin{enumerate}[(i)]
    \item any object in $\bfC$ is compact in $\bfE$, and
    \item any object in $\bfE$ is a filtered colimit of objects in $\bfC$,
\end{enumerate}
then the canonical extension $\Ind(\bfC) \to \bfE$, coming from the universal property of $\Ind(\bfC)$, is an equivalence of categories.
\end{lemma}

To avoid size issues, we assume from now on that $\bfC$ is an (essentially) small category. The fact that the presheaf category $\Set^{\bfC^{op}}$ is the free cocompletion of $\bfC$ leads to an alternative description of $\Ind(\bfC)$ that is sometimes easier to work with. Namely, we can think of $\Ind(\bfC)$ as the full subcategory of $\Set^{\bfC^{op}}$ consisting of those presheaves which are filtered colimits of representables. If $\bfC$ is small and has finite colimits, as will be the case in all of our examples, then these are exactly the functors $\bfC^{op} \to \Set$ that send the finite colimits of $\bfC$ to limits in $\Set$ (see Théorème 8.3.3.(v) of \cite[Exposé 1]{grothendieck1972theorie}), i.e.
\[ \Ind(\bfC) \simeq \mathrm{lex}(\bfC^{op},\Set), \]
where the right-hand side stands for the category of left exact functors. From this description, one sees immediately that $\Ind(\bfC)$ has all small limits and that the inclusion $\Ind(\bfC) \to \Set^{\bfC^{op}}$ preserves these. The category $\Ind(\bfC)$ also has all colimits in this case. Namely, finite coproducts and pushouts can be computed ``levelwise'' in $\bfC$ as described in \cite[Appendix 4]{ArtinMazur1986Etale}, while filtered colimits exist as mentioned above. Note however, that the inclusion $\Ind(\bfC) \to \Set^{\bfC^{op}}$ does not preserve all colimits, but only filtered ones.

One sees dually that if $\bfC$ is small and has all finite limits, then 
\[ \Pro(\bfC) \simeq \mathrm{lex}(\bfC,\Set)^{op}. \]
As above, it follows that $\Pro(\bfC)$ is complete and cocomplete in this case.

Another consequence of the fact that finite coproducts and pushouts in $\Ind(\bfC)$ are computed ``levelwise'' is the following: if $F \colon \bfC \to \bfE$, with $\bfE$ cocomplete, preserves finite colimits, then its extension $\wt F \colon \Ind(\bfC) \to \bfE$ given by the universal property also preserves finite colimits. Since it also preserves filtered colimits, we conclude that it preserves all colimits. In fact, more is true. The above description of $\Ind(\bfC)$ as $\mathrm{lex}(\bfC^{op}, \Set)$ allows us to construct a right adjoint $R$ of $\wt F$. Namely, if we define $R(E)(c) := \Hom(Fc,E)$, then $R(E) \colon \bfC^{op} \to \Set$ is left exact, hence $R$ defines a functor $\bfE \to \Ind(\bfC)$. Adjointness follows from the Yoneda lemma. We therefore see that, up to unique natural isomorphism, there is a 1-1 correspondence between finite colimit perserving functors $\bfC \to \bfE$ and functors $\Ind(\bfC) \to \bfE$ that have a right adjoint.

There are two important examples of adjunctions obtained in this way that we would like to mention here. The first one is the ind-completion functor. If $\bfE$ is a cocomplete category and $\bfC$ a full subcategory closed under finite colimits, then the inclusion $\bfC \subset \bfE$ induces an adjunction
\[U : \Ind(\bfC) \rightleftarrows \bfE : \wh{(\cdot)}_{\Ind}\]
whose right adjoint we call \emph{ind-completion (relative to $\bfC$)} or \emph{ind-$\bfC$ completion}. Dually, if $\bfE$ is complete and $\bfC$ is a full subcategory closed under finite limits, then we obtain an adjunction
\[\wh{(\cdot)}_{\Pro} : \bfE \rightleftarrows \Pro(\bfC) : U \]
whose left adjoint we call \emph{pro-completion (relative to $\bfC$)} or \emph{pro-$\bfC$ completion}. In many examples, $\bfC$ is the full subcategory of $\bfE$ consisting of objects that are ``finite'' in some sense, and this left adjoint is better known as the profinite completion functor. For instance, in the case of groups, this functor $\wh{(\cdot)}_{\Pro} \colon \Grp \to \Pro(\Fin\Grp)$ is the well-known profinite completion functor for groups.

The other important example is about cotensors in ind-categories. Suppose $\bfC$ is a small simplicial category that has all finite colimits and tensors with finite simplicial sets, and that furthermore these tensors commute with these finite colimits. We will call $\bfC$ \emph{finitely tensored} if this is the case (cf.~\autoref{definition:FinitelyTensored} for a precise definition). If $X$ is a simplicial set, then we can write it as $\colim_i X_i$ where $i$ ranges over all finite simplicial subsets $X_i \subset X$. Define $- \otimes X \colon \bfC \to \Ind(\bfC)$ by $c \otimes X = \{c \otimes X_i\}_i$. This functor preserves finite colimits since these are computed ``levelwise'' in $\Ind(\bfC)$, hence it extends to a functor $- \otimes X \colon \Ind(\bfC) \to \Ind(\bfC)$ that has a right adjoint $(-)^X$. These define tensors and cotensors by arbitrary simplicial sets on $\Ind(\bfC)$. In particular, $\Ind(\bfC)$ is a simplicial category that is complete, cocomplete, tensored and cotensored (note the similarity with Proposition 4.10 of \cite{BarneaSchlank2016Weak}). The dual of this statement says that for any small simplicial category $\bfC$ that has finite limits and cotensors with finite simplicial sets, and in which these finite cotensors commute with finite limits in $\bfC$, the pro-category $\Pro(\bfC)$ is a simplicial category that is complete, cocomplete, tensored and cotensored. We call $\bfC$ \emph{finitely cotensored} in this case.

Let us return to the basic definition \eqref{equation:DefinitionHomSetIndCategory} of morphisms in $\Ind(\bfC)$. If $C = \{c_i\}$ and $D = \{d_i\}$ are objects indexed by the same filtered category $I$, then any natural transformation with components $f_i \colon c_i \to d_i$ represents a morphism in $\Ind(\bfC)$. Morphisms of this type (or more precisely, morphisms represented in this way) will be called \emph{level maps} or \emph{strict maps}. Up to isomorphism, any morphism in $\Ind(\bfC)$ has such a strict representation (see Corollary 3.2 of \cite[Appendix]{ArtinMazur1986Etale}). One can define the notion of a ``level'' diagram or ``strict'' diagram in a similar way. Given an indexing category $K$, a conceptual way of thinking about these is through the canonical functor
\[ \Ind(\bfC^K) \to \Ind(\bfC)^K. \]
A strict diagram can be thought of as an object in the image of this functor. If $K$ is a finite category and $\bfC$ has all finite colimits, then the above functor is an equivalence of categories (\cite[\S 4]{Meyer1980Approximation}). This shows in particular that, up to isomorphism, any finite diagram in $\Ind(\bfC)$ is a strict diagram if $\bfC$ is small and has finite colimits.

In our context, the following extension of Meyer's result is important. Suppose that $K$ is a category which can be written as a union of a sequence of finite full subcategories
\[ K_0 \subset K_1 \subset K_2 \subset \ldots \subset K = \cup_{n \in \bbN} K_n. \]
Let $\bfC$ be a small category that has finite colimits. Then any functor $f \colon K_n \to \bfC$ has a left Kan extension $g \colon K \to \bfC$ defined in terms of finite colimits as in (the dual to) Theorem X.3.1 of \cite{MacLane1971CategoriesWorking}. For $X \colon K \to \bfC$, write $\sk_n X$ for the left Kan extension of the restriction of $X$ to $K_n$. We call $X$ \emph{$n$-skeletal} if the canonical map $\sk_n X \to X$ is an isomorphism, and \emph{skeletal} if this is the case for some $n$. The full subcategory $\sk(\bfC^K) \subset \bfC^K$ spanned by the skeletal functors $K \to C$ can be viewed as a full subcategory of $\Ind(\bfC)^K$ via the inclusion $\bfC \hookrightarrow \Ind(\bfC)$. Note that for any $X$ in $\Ind(\bfC)^K$, we have $X = \colim_n \sk_n X$. Exactly as in (the dual of) the proof of Proposition 7.4.1 of \cite{BarneaHarpazHorel2017}, the result of \cite[\S 4]{Meyer1980Approximation} mentioned above can be used to show that the hypotheses of the recognition principle for ind-categories are satisfied, hence that the induced functor $\Ind(\sk(\bfC^K)) \to \Ind(\bfC)^K$ is an equivalence of categories. In fact, the assumption that $K$ is a union of a \emph{sequence} of finite full subcategories is irrelevant, and the following more general result, which we write down for future reference, can be proved by the same argument. Note that a category $K$ can be written as a union of finite full subcategories if and only if for any $k, k' \in K$, the set $\Hom_K(k,k')$ is finite.

\begin{theorem}\label{thm:IndValuedPresheavesIsIndSkeletalPresheaves}
Let $\bfC$ be a small category that has finite colimits, and let $K$ be a small category that can be written as a union of finite full subcategories. Write $\sk(\bfC^K)$ for the full subcategory of $\bfC^K$ of those functors $K \to \bfC$ that are isomorphic to the left Kan extension of a functor $K' \to \bfC$ for some finite full subcategory $K' \subset K$. Then $\Ind(\sk(\bfC^K)) \simeq \Ind(\bfC)^K$.
\end{theorem}

This theorem recovers the well-known equivalence $\Ind(\FinsSet) \simeq \s\Set$ when applied to $\Delta^{op} = \cup_n \Delta^{op}_{\leq n}$ and $\bfC = \Fin\Set$. Note that we already (implicitly) used this equivalence when we defined tensors by simplicial sets for ind-categories above.

We can also apply the dual of this theorem to the same categories $K = \Delta^{op}$ and $\bfC = \Fin\Set$. Write $\wh \Set = \Pro(\Fin\Set)$ for the category of profinite sets, which is well known to be equivalent to the category of Stone spaces $\Stone$. Since we want to apply the dual of \autoref{thm:IndValuedPresheavesIsIndSkeletalPresheaves}, we need to work with right Kan extensions instead of left Kan extensions. In particular, we obtain the full subcategory of $\Fin\Set^{\Delta^{op}}$ on those simplicial sets that are the right Kan extension of some functor $\Delta^{op}_{\leq n} \to \Fin\Set$. These are exactly the coskeletal degreewise finite simplicial sets, i.e.~the lean simplicial sets. In particular, the theorem above recovers the equivalence $\Pro(\LsSet) \simeq \s\wh\Set$ proved in Proposition 7.4.1 of \cite{BarneaHarpazHorel2017}.

An example that plays an important role in \autoref{sec:CompleteSegalProfiniteVSProfiniteQuasiCats} is that of bisimplicial (profinite) sets. The dual of the above theorem shows that the category of bisimplicial profinite sets $\biswhSet = \s\wh\Set^{\Delta^{op}} \cong \wh\Set^{\Delta^{op} \times \Delta^{op}}$ is canonically equivalent to the category $\Pro(\LsSet^{(2)})$ for a certain full subcategory $\LsSet^{(2)}$ of the category of bisimplicial sets $\bisSet = \s\Set^{\Delta^{op}}$. This category $\LsSet^{(2)}$ consists of those bisimplicial sets that are isomorphic to the right Kan extension of a functor $\Delta^{op}_{\leq t} \times \Delta^{op}_{\leq n} \to \Fin\Set$ along the inclusion $\Delta^{op}_{\leq t} \times \Delta^{op}_{\leq n} \hookrightarrow \Delta^{op} \times \Delta^{op}$ for some $t,n \in \bbN$. We will refer to such bisimplicial sets as \emph{doubly lean}.


\section{The completed model structure on \texorpdfstring{$\Ind(\bfC)$}{Ind(C)}}\label{sec:IndCompletedModelStructure}

In this section, we will describe our construction of the model structure on $\Ind(\bfC)$, where $\bfC$ is what we call a ``cofibration test category''. In \autoref{sec:ProCompletedModelStructure}, we will dualize this construction to the context of pro-categories. After that, we will study the functorial behaviour of the construction in \autoref{sec:QuillenPairs} and discuss Bousfield localizations in \autoref{sec:BousfieldLocalizations}.

Throughout these sections, the terms ``weak equivalence'' and ``fibration'' of simplicial sets refer to either the classical Kan-Quillen model structure or to the Joyal model structure. When we say that a model category is simplicial, then this can either mean that it is enriched over the Kan-Quillen model structure or over the Joyal model structure.

We wish to single out the definition of being finitely tensored, since it occurs many times throughout this paper.

\begin{definition}\label{definition:FinitelyTensored}
	Let $\bfC$ be simplicial category. Then $\bfC$ is called \emph{finitely tensored} if
	\begin{enumerate}[(i)]
		\item\label{item1:FinitelyTensored} it admits finite colimits,
		\item\label{item2:FinitelyTensored} it admits tensors by finite simplicial sets, and
		\item\label{item3:FinitelyTensored} these commute with each other, meaning that the canonical map
		\[ \colim_i (c_i \otimes X) \to (\colim_i c_i) \otimes X\]
		is an isomorphism for any finite diagram $\{c_i\}$ in $\bfC$ and any finite simplicial set $X$.
	\end{enumerate}
\end{definition}

\begin{remark}
	It is worth pointing out that condition \eqref{item3:FinitelyTensored} is equivalent to asking that the finite colimits of \eqref{item1:FinitelyTensored} are \emph{enriched colimits}; that is, for any finite diagram $\{c_i\}$ in $\bfC$ and any object $d$ in $\bfC$, the canonical map $\colim_i \Map(c_i, d) \to \Map(\colim_i c_i, d)$ is an isomorphism of simplicial sets.
\end{remark}

As explained in \autoref{ssec:ProAndIndCategories}, if $\bfC$ is finitely tensored, then the category $\Ind(\bfC)$ is a tensored and cotensored simplicial category that is both complete and cocomplete. We will endow $\bfC$ with some additional structure, that of a ``cofibration test category'', and show that it induces a simplicial model structure on $\Ind(\bfC)$ in \autoref{theorem:ModelStructureIndC} below.

\begin{definition} \label{definition:CofTestCat}
A \emph{cofibration test category} $(\bfC, \bfT)$ consists of a small finitely tensored simplicial category $\bfC$, a full subcategory $\bfT$ of \emph{test objects} and two classes of maps in $\bfT$ called \emph{cofibrations}, denoted $\cofarrow$, and \emph{trivial cofibrations}, denoted $\trivcofarrow$, both containing all isomorphisms, that satisfy the following properties:
\begin{enumerate}[(1)]
    \item\label{CofTestCat:TestIsCofibrant} The initial object $\varnothing$ is a test object, and for every test object $t \in \bfT$, the map $\varnothing \to t$ is a cofibration.
    \item\label{CofTestCat:TestObjectsClosedUnderTensors} For every cofibration between test objects $s \cofarrow t$ and cofibration between finite simplicial set $U \cofarrow V$, the pushout-product map $t \otimes U \cup_{s \otimes U} s \otimes V \to t \otimes V$ is a cofibration between test objects which is trivial if either $s \cofarrow t$ or $U \cofarrow V$ is so.
    \item\label{CofTestCat:TestingTrivialCofibration} A morphism $r \to s$ in $\bfT$ is a trivial cofibration if and only if it is a cofibration and $\Map(t,r) \to \Map(t,s)$ is a weak equivalence of simplicial sets for every $t \in \bfT$.
    \item\label{CofTestCat:FibrantMappingSpaces} Any object $c \in \bfC$ has the right lifting property with respect to trivial cofibrations.
\end{enumerate}
\end{definition}

\begin{remark}\label{remark:FibrantMappingSpacesIndC}
Note that property \eqref{CofTestCat:FibrantMappingSpaces} implies that $\Map(t,C)$ is fibrant for every $t \in \bfT$ and $C \in \Ind(\bfC)$. Namely, writing $C$ as a filtered colimit $\colim_i c_i$ with $c_i \in \bfC$ for every $i$, we see that $\Map(t,C) = \colim_i \Map(t,c_i)$. Hence it suffices to show that $\Map(t,c)$ is fibrant for every object $c$ in $\bfC$. This is equivalent to $c$ having the right lifting property with respect to certain maps of the form $t \otimes \Lambda^n_k \to t \otimes \Delta^n$, which is indeed the case by items \eqref{CofTestCat:TestIsCofibrant}, \eqref{CofTestCat:TestObjectsClosedUnderTensors} and \eqref{CofTestCat:FibrantMappingSpaces}.
\end{remark}

\begin{remark}\label{remark:CofTestCatJoyalOrKQ}
The definition of a cofibration test category depends on whether we work with the Kan-Quillen model structure $\s\Set_{KQ}$ or the Joyal model structure $\s\Set_J$. However, since $\s\Set_{KQ}$ is a left Bousfield localization of $\s\Set_J$, any cofibration test category with respect to $\s\Set_{KQ}$ is also a cofibration test category with respect to $\s\Set_J$. To see this, suppose that $(\bfC,\bfT)$ is a cofibration test category with respect to $\s\Set_{KQ}$. It is clear that items \eqref{CofTestCat:TestIsCofibrant}, \eqref{CofTestCat:TestObjectsClosedUnderTensors} and \eqref{CofTestCat:FibrantMappingSpaces} also hold with respect to $\s\Set_J$. For item \eqref{CofTestCat:TestingTrivialCofibration}, note that the map $\Map(t,r) \to \Map(t,s)$ is a map between Kan complexes by \autoref{remark:FibrantMappingSpacesIndC}, hence that it is a weak equivalence in $\s\Set_J$ if and only if it is in $\s\Set_{KQ}$.
\end{remark}

We will often write $\bfC$ for a cofibration test category $(\bfC,\bfT)$, omitting the full subcategory of test objects $\bfT$ from the notation. We will write $\cof(\bfC)$ for the set of cofibrations. Note that this is a subset of the morphisms of $\bfT$.

The role of the test objects $t \in \bfT$ is to detect the weak equivalences in $\Ind(\bfC)$ ``from the left''. More precisely, the weak equivalences in $\Ind(\bfC)$ will be those arrows $C \to D$ for which $\Map(t,C) \to \Map(t,D)$ is a weak equivalence for every $t \in \bfT$. For this reason, we will call an arrow $c \to d$ in $\bfC$ for which $\Map(t,c) \to \Map(t,d)$ is a weak equivalence for every $t \in \bfT$ a \emph{weak equivalence}, and denote such arrows by $\wearrow$. We write $\we(\bfC)$ for the set of weak equivalences in $\bfC$. Using this terminology, item \eqref{CofTestCat:TestingTrivialCofibration} of the above definition can be rephrased as saying that the trivial cofibrations are precisely the cofibrations that are weak equivalences. In particular, the set of trivial cofibrations in a cofibration test category $\bfC$ is fully determined by the full subcategory $\bfT$ and the set $\cof(\bfC)$.

Let us look at a few examples. Note that we will discuss more interesting examples in \autoref{sec:ProCompletedModelStructure}, where we consider fibration test categories, the dual of cofibration test categories.

\begin{example}\label{example:InheritedCofibrationTestCategory}
Suppose $\bfE$ is a simplicial model category in which every object is fibrant and let $\bfC \subset \bfE$ be a (small) full subcategory closed under finite colimits and finite tensors. If we define $\bfT$ to be the full subcategory on the cofibrant objects, then $(\bfC,\bfT)$ forms a cofibration test category where the (trivial) cofibrations are the (trivial) cofibrations of $\bfE$ between objects of $\bfT$. We say that $\bfC$ \emph{inherits} this structure of a cofibration test category from $\bfE$. Properties \eqref{CofTestCat:TestIsCofibrant}, \eqref{CofTestCat:TestObjectsClosedUnderTensors} and \eqref{CofTestCat:FibrantMappingSpaces} of \autoref{definition:CofTestCat} follow directly from the fact that $\bfE$ is a (simplicial) model category and the fact that any object in $\bfE$ is fibrant. For one direction of property \eqref{CofTestCat:TestingTrivialCofibration}, note that since all object in $\bfE$ are fibrant, the functor $\Map(t,-)$ preserves weak equivalences for any cofibrant object $t$. For the converse direction, note that a cofibration $r \cofarrow s$ is trivial if and only if it is mapped to an isomorphism in the homotopy category $\Ho(\bfE)$. By the Yoneda lemma applied to the full subcategory $\Ho(\bfT) \subset \Ho(\bfE)$ spanned by the objects of $\bfT$, this is equivalent to $\Hom_{\Ho(\bfE)}(t,r) \to \Hom_{\Ho(\bfE)}(t,s)$ being a weak equivalence for every $t$. Since $\Map(t,r) \to \Map(t,s)$ is a weak equivalence by assumption and $\Hom_{\Ho(\bfE)}(t,\mhyphen)$ equals the set of path components of (the maximal Kan complex contained in) $\Map(t,\mhyphen)$, this is indeed the case.
\end{example}

\begin{example}\label{example:RestrictingToLessTestObjects}
Suppose that a cofibration test category $(\bfC, \bfT)$ is given, and let $\bfT' \subset \bfT$ be a full subcategory such that $\varnothing \in \bfT'$ and such that for any cofibration $s \cofarrow t$ between objects of $\bfT'$ and any cofibration $U \cofarrow V$ in $\FinsSet$, the object $t \otimes U \cup_{s \otimes U} s \otimes V$ is again in $\bfT'$. We will call such a full subcategory $\bfT' \subset \bfT$ \emph{closed under finite pushout-products}. Then $(\bfC, \bfT')$ is again a cofibration test category if we define the (trivial) cofibrations to be those of $(\bfC,\bfT)$ between objects of $\bfT'$. All items of \autoref{definition:CofTestCat} are straightforward to show except possibly property \eqref{CofTestCat:TestingTrivialCofibration}. The ``only if'' direction follows immediately. For the ``if'' direction of \eqref{CofTestCat:TestingTrivialCofibration}, suppose $r \cofarrow s$ is a map in $\bfT'$ that is a cofibration with the property that $\Map(t,r) \to \Map(t,s)$ is a weak equivalence for any $t \in \bfT'$. Applying this to $t=r$ and $t=s$ and using that these mapping spaces are fibrant, we obtain left and right homotopy inverses of $r \cofarrow s$, where homotopies in $\bfT'$ are defined using the tensor $- \otimes \Delta^1$ (in the case of $\s\Set_{KQ}$) or $- \otimes H$ (in the case $\s\Set_J$, where $H$ is as in the proof of \autoref{lem:SkeletalJ}). Since $\Map(t,\mhyphen)$ is a simplicial functor it preserves these homotopies, showing that $\Map(t,r) \to \Map(t,s)$ is homotopy equivalence for every $t \in \bfT$. We conclude that $r \cofarrow s$ is a trivial cofibration in $\bfT$ and hence a trivial cofibration in $\bfT'$ by definition.
\end{example}

\begin{example}\label{example:ConvenientCofTestCatSpaces}
Let $\Top$ be a convenient category of topological spaces, such as $k$-spaces or compactly generated (weak) Hausdorff spaces. The Quillen model structure on $\Top$ is a simplicial model structure, in which tensors are given by $C \otimes X = C \times |X|$ for any $C \in \Top$ and $X \in \s\Set$. Let $\bfC \subset \Top$ be any small full subcategory of $\Top$ that is closed under finite colimits and finite tensors, and moreover contains the space $|X|$ for any finite simplicial set $X$. Define $\bfT \subset \bfC$ to be the full subcategory consisting of the objects $|X|$ for any finite simplicial set $X$, and define a map $|X| \to |Y|$ in $\bfT$ to be a (trivial) cofibration if it is the geometric realization of a (trivial) cofibration $X \cofarrow Y$ in the Kan-Quillen model structure on $\s\Set$. Using that there are natural isomorphisms $|Y| \otimes V \cong |Y \times V|$ and $|X| \otimes V \cup_{|X| \otimes U} |Y| \otimes U \cong |X \times V \cup_{X \times U} Y \times V|$ for any pair of maps $X \to Y$ and $U \to V$ in $\s\Set$, it is straightforward to verify that $(\bfC,\bfT)$ is a cofibration test category in the sense of \autoref{definition:CofTestCat} (with respect to $\s\Set_{KQ}$). This example will be studied further in \autoref{sec:ConvenientCategorySpaces}.
\end{example}

For a cofibration test category $\bfC$, we will write $I$ for the image of the set of cofibrations of $\bfC$ in $\Ind(\bfC)$, and $J$ for the image of the set of trivial cofibrations of $\bfC$ in $\Ind(\bfC)$. Identifying $\bfC$ with its image in $\Ind(\bfC)$, we can write
\[I = \{f \colon s \to t \mid f \text{ is a cofibration in } \bfC \} = \cof(\bfC) \]
and
\[J = \{f \colon s \to t \mid f \text{ is a trivial cofibration in } \bfC \} = \cof(\bfC) \cap \we(\bfC). \]
Recall that the sets of (trivial) cofibrations $\cof(\bfC)$ and $\cof(\bfC) \cap \we(\bfC)$ in $(\bfC,\bfT)$ are both contained in $\bfT$; that is, any (trivial) cofibration is a map between test objects.

The sets $I$ and $J$ are generating (trivial) cofibrations for a model structure on $\Ind(\bfC)$ in which the weak equivalences are as above.

\begin{theorem}\label{theorem:ModelStructureIndC}
Let $\bfC$ be a cofibration test category. Then $\Ind(\bfC)$ carries a cofibrantly generated (hence combinatorial) simplicial model structure, the \emph{completed model structure}, where a map $C \to D$ is a weak equivalence if and only if $\Map(t,C) \to \Map(t,D)$ is a weak equivalence for every $t \in \bfT$. A set of generating cofibrations (generating trivial cofibrations) is given by $I$ ($J$, respectively). Every object is fibrant in this model structure.
\end{theorem}

\begin{remark}
As mentioned in \autoref{remark:CofTestCatJoyalOrKQ}, the definition of a cofibration test category depends on whether we work with the Joyal model structure or the Kan-Quillen model structure on $\s\Set$. In the first case, the model structure on $\Ind(\bfC)$ will be $\s\Set_J$-enriched, while in the latter case, it will be $\s\Set_{KQ}$-enriched.
\end{remark}

The proof uses the following lemmas.

\begin{lemma}\label{lemma:WeakEquivalencesStableUnderFilCol}
Let $\bfC$ be a cofibration test category. The weak equivalences of $\Ind(\bfC)$ as defined in \autoref{theorem:ModelStructureIndC} are stable under filtered colimits.
\end{lemma}

\begin{proof}
Let $\{C_i \wearrow D_i\}$ be a levelwise weak equivalence between filtered diagrams in $\Ind(\bfC)$ and let $t \in \bfT$. Then $\Map(t,\colim_i C_i) \to \Map(t,\colim_i D_i)$ is the filtered colimit of the maps $\Map(t,C_i) \to \Map(t,D_i)$ since $t$ is compact in $\Ind(\bfC)$, which are weak equivalences by assumption. The proof therefore reduces to the statement in $\s\Set$ that a filtered colimit of weak equivalences, indexed by some filtered category $I$, is again a weak equivalence. This can be proved for the Kan-Quillen and Joyal model structure in exactly the same way. Namely, this is equivalent to the statement that the functor $\colim \colon \s\Set^I \to \s\Set$, where $\s\Set^I$ is endowed with the projective model structure, preserves weak equivalences. To see that this is the case, factor $\{X_i\} \wearrow \{Y_i\}$ in $\s\Set^I$ as a projective trivial cofibration $\{X_i\} \trivcofarrow \{Z_i\}$ followed by a pointwise trivial fibration $\{Z_i\} \trivfibarrow \{Y_i\}$. Then $\colim X_i \to \colim Z_i$ is again a trivial cofibration, so in particular a weak equivalence. Furthermore, since the generating cofibrations $\partial \Delta^n \to \Delta^n$ in $\s\Set$ are maps between compact objects, we see that $\colim Z_i \to \colim Y_i$ must have the right lifting property with respect to these maps, i.e.~it is a trivial fibration. We conclude that $\colim \colon \s\Set^I \to \s\Set$ preserves weak equivalences.
\end{proof}

\begin{lemma}\label{lemma:PullbackMappingPropertyForIndC}
Let $\bfC$ be a cofibration test category, let $s \cofarrow t$ be a cofibration in $\bfC$, i.e.~a map in $I$, and let $C \to D$ be an arrow in $\Ind(\bfC)$ which has the right lifting property with respect to all maps in $J$. Then $\Map(t,C) \to \Map(s,C) \times_{\Map(s,D)} \Map(t,D)$ is a fibration, which is trivial if either $s \cofarrow t$ is trivial or if $C \to D$ is a weak equivalence in the sense of \autoref{theorem:ModelStructureIndC}.
\end{lemma}

\begin{proof}
Let $M$ be a set of trivial cofibrations in $\FinsSet$ such that a map between fibrant objects in $\s\Set$ is a fibration if and only if it has the right lifting property with respect to the maps in $M$. For the Kan-Quillen model structure, one can take the set of horn inclusions, while for $\s\Set_J$, the set $M$ from \autoref{lem:SkeletalJ} works. By \autoref{remark:FibrantMappingSpacesIndC}, for any test object $t \in \bfT$ and any $C \in \Ind(\bfC)$ the simplicial set $\Map(t,C)$ is fibrant. For any $t \in \bfT$ and any map $U \trivcofarrow V$ in $M$, the map $t \otimes U \to t \otimes V$ is in $J$ by items \eqref{CofTestCat:TestIsCofibrant} and \eqref{CofTestCat:TestObjectsClosedUnderTensors} of \autoref{definition:CofTestCat}. By adjunction, we conclude that for any $C \to D$ that has the right lifting property with respect to maps in $J$, the map $\Map(t,C) \to \Map(t,D)$ is a fibration. If we are given a map $s \cofarrow t$ in $I$, then $\Map(s,C) \times_{\Map(s,D)} \Map(t,D)$ is fibrant because the map to $\Map(t,D)$ is the pullback of the fibration $\Map(s,C) \fibarrow \Map(s,D)$. By a similar argument as above, $\Map(t,C) \to \Map(s,C) \times_{\Map(s,D)} \Map(t,D)$ is a fibration. The same argument with the set of boundary inclusions $\{ \partial \Delta^n \to \Delta^n\}$ instead of $M$ shows that $\Map(t,C) \to \Map(s,C) \times_{\Map(s,D)} \Map(t,D)$ is a trivial fibration if $s \cofarrow t$ is in $J$. If $C \to D$ is a weak equivalence, then the maps $\Map(s,C) \to \Map(s,D)$ and $\Map(t,C) \to \Map(t,D)$ are weak equivalences by definition, hence trivial fibrations by the above. As indicated in the diagram
\[\begin{tikzcd}[column sep = scriptsize]
\Map(t,C) \arrow[drr, bend left = 13] \arrow[ddr, bend right, "\sim" rot320, two heads] \arrow[dr, two heads] & & \\
& \Map(s,C) \times_{\Map(s,D)} \Map(t,D) \arrow[dr,phantom, very near start, "\lrcorner"] \arrow[r] \arrow[d, two heads] & \Map(s,C) \arrow[d, "\sim" rot90, two heads] \\
& \Map(t,D) \arrow[r] & \Map(s,D).
\end{tikzcd}\]
the map $\Map(t,C) \to \Map(s,C) \times_{\Map(s,D)} \Map(t,D)$ is a trivial fibration by the 2 out of 3 property.
\end{proof}

\begin{lemma}\label{lemma:TrivCofibStableUnderPushout}
Let $\bfC$ be a cofibration test category and let $s \trivcofarrow t$ be a trivial cofibration in $\bfC$. Then any pushout of $s \trivcofarrow t$ in $\Ind(\bfC)$ is a weak equivalence in the sense of \autoref{theorem:ModelStructureIndC}.
\end{lemma}

\begin{proof}
The following proof works if $\bfC$ is a cofibration test category with respect to the Kan-Quillen model structure on $\s\Set$. The same proof works in the case that $\bfC$ is a cofibration test category with respect to $\s\Set_J$ if one replaces every instance of $\Delta^1$ by $H$, where $H$ is as in the proof of \autoref{lem:SkeletalJ}.

We will first show that $i \colon s \trivcofarrow t$ is a deformation retract. By item \eqref{CofTestCat:FibrantMappingSpaces} of \autoref{definition:CofTestCat}, there exists a lift in
\[\begin{tikzcd}
s \ar[r,"="] \ar[d,tail, "\sim" rot90,"i"] & s \\
t, \ar[ur,"r"',dashed]
\end{tikzcd}\]
i.e.~a retract $r$ of $i$. By \eqref{CofTestCat:TestObjectsClosedUnderTensors} and \eqref{CofTestCat:FibrantMappingSpaces} of \autoref{definition:CofTestCat}, there exists a lift $F$ in
\[\begin{tikzcd}[column sep = huge]
s \otimes \Delta^1 \cup_{s \otimes \partial \Delta^1} t \otimes \partial \Delta^1 \ar[r,"i \cup (\id_t{,}ir)"] \ar[d,tail, "\sim" rot90] & t \\
t \otimes \Delta^1 \ar[ur,dashed,"F"'] & 
\end{tikzcd}\]
This lift $F \colon t \otimes \Delta^1 \to t$ is a deformation retract.

Now let a pushout square
\begin{equation}\label{diagram:PushoutDeformRetract}\begin{tikzcd}
s \ar[d,tail, "\sim" rot90,"i"]\ar[r,"f"] & C \ar[d,"j"]  \\
t \ar[r,"g"] & D \ar[ul,phantom, very near start, "\ulcorner"]
\end{tikzcd}\end{equation}
be given, where $s \to C$ is any map in $\Ind(\bfC)$. The maps $fr \colon t \to C$ and $\id_C \colon C \to C$ give, by the universal property of the pushout, a retract $r'$ of $j \colon C \to D$. Since tensors preserve colimits, we see that $D \otimes \Delta^1$ is the pushout of $t \otimes \Delta^1$ and $C \otimes \Delta^1$ along $s \otimes \Delta^1$. Then $g \circ F \colon t \otimes \Delta^1 \to D$ and $C \otimes \Delta^1 \to C \otimes * \cong C \xrightarrow{j} D$ give, by the universal property of the pushout, a map $G \colon D \otimes \Delta^1 \to D$. Write $\iota_0,\iota_1 \colon D \to D \otimes \Delta^1$ for the endpoint inclusions. It follows from the universal property of the pushout \eqref{diagram:PushoutDeformRetract} that $\iota_0 G = \id_D$ while $\iota_1 G = jr'$, i.e.~$G$ is a deformation retract.

Now let $u \in \bfT$ be any test object. We deduce from the existence of the deformation retract $G$ that $\Map(u,C) \to \Map(u,D)$ is the inclusion of a deformation retract, hence a weak equivalence.
\end{proof}

\begin{proof}[Proof of \autoref{theorem:ModelStructureIndC}]
We check all the four assumptions of Kan's recognition theorem as spelled out in \cite[Theorem 11.3.1]{Hirschhorn2003Model}. The weak equivalences satisfy the two out of three property and are closed under retracts since this holds for the weak equivalences in $\s\Set$.

\par (1) Since all objects of $\bfC$ are compact in $\Ind(\bfC)$, the sets $I$ and $J$ permit the small object argument.

\par (2) It suffices to prove that any transfinite composition of pushouts of maps in $J$ is a weak equivalence. This follows immediately from  \autoref{lemma:WeakEquivalencesStableUnderFilCol} and \autoref{lemma:TrivCofibStableUnderPushout}.

\par (3) We need to show that any map having the right lifting property with respect to maps in $I$ has the right lifting property with respect to maps in $J$ and is a weak equivalence. The first of these follows since $J \subset I$. To see that any map that has the right lifting property with respect to maps in $I$ is a weak equivalence, let such a map $C \to D$ be given. Note that $t \otimes \partial \Delta[n] \to t \otimes \Delta[n]$ is in $I$ for any $t \in \bfT$ and $n \geq 0$ by items \eqref{CofTestCat:TestIsCofibrant} and \eqref{CofTestCat:TestObjectsClosedUnderTensors} of \autoref{definition:CofTestCat}. This implies that $\Map(t,C) \to \Map(t,D)$ is a trivial fibration for any $t \in \bfT$ and in particular that $C \to D$ is a weak equivalence.

\par (4) We need to show that if $C \to D$ has the right lifting property with respect to maps in $J$ and is a weak equivalence, then it has the right lifting property with respect to maps in $I$. Let $s \cofarrow t$ in $I$ be given. Then $\Map(t,C) \trivfibarrow \Map(s,C) \times_{\Map(s,D)} \Map(t,D)$ is a trivial fibration by \autoref{lemma:PullbackMappingPropertyForIndC}, and in particular surjective on $0$-simplices. In particular, $C \to D$ has the right lifting property with respect to $s \cofarrow t$.

The fact that this model structure is simplicial follows from \autoref{lemma:PullbackMappingPropertyForIndC}. By \eqref{CofTestCat:FibrantMappingSpaces}, all objects in $\bfC \subset \Ind(\bfC)$ are fibrant. Since the generating trivial cofibrations are maps between compact objects and any $C \in \Ind(\bfC)$ is a filtered colimit of objects in $\bfC$, it follows that all objects in $\Ind(\bfC)$ are fibrant.
\end{proof}

\begin{example}\label{example:ConvenientModelCatOfSpaces}
Let $(\bfC,\bfT)$ be the cofibration test category from \autoref{example:ConvenientCofTestCatSpaces}. The model structure on $\Ind(\bfC)$ obtained by applying \autoref{theorem:ModelStructureIndC} turns out to be Quillen equivalent to the Kan-Quillen model structure on $\s\Set$ and the Quillen model structure on $\Top$. More precisely, there is a canonical way to factor the geometric realization functor for simplicial sets $|\cdot| \colon \s\Set \to \Top$  as a composite $\s\Set \to \Ind(\bfC) \to \Top$, where both of these functors are left Quillen equivalences. This will be proved in \autoref{proposition:ConventientCategoryTopologicalEquivalences}.
\end{example}

\begin{example}
For this example, $\Top$ is again a convenient category of spaces as in \autoref{example:ConvenientCofTestCatSpaces}. If $\mathcal{P}$ is a topological operad, then the category $\mathcal{P}\mhyphen\mathbf{Alg}$ of $\mathcal{P}$-algebras admits a model structure, obtained through transfer along the free-forgetful adjunction $F : \Top \rightleftarrows \mathcal{P} \mhyphen \mathbf{Alg} : U$. In particular, any object is fibrant in this model structure. This category is $\Top$-enriched, since one can view $\Hom_{\mathcal{P} \mhyphen \mathbf{Alg}}(S,T)$ as a subspace of $\ul\Hom_{\Top}(US,UT)$. For any topological space $X$ and any $\mathcal{P}$-algebra $S$, one can endow the space $S^X$ with the ``pointwise'' structure of a $\mathcal{P}$-algebra. By restricting the usual homeomorphism coming from the cartesian closed structure on $\Top$, we obtain a natural homeomorphism $\ul \Hom_{\mathcal{P} \mhyphen \mathbf{Alg}}(S,T^X) \cong \ul \Hom_{\Top}(X, \ul\Hom_{\mathcal{P} \mhyphen \mathbf{Alg}}(S,T))$. One can furthermore show that $-^X \colon \mathcal{P} \mhyphen \mathbf{Alg} \to \mathcal{P} \mhyphen \mathbf{Alg}$ has a left adjoint that makes $\mathcal{P} \mhyphen \mathbf{Alg}$ into a tensored and cotensored topological category. In particular, it can be viewed as a tensored and cotensored simplicial category. Since the cotensors, fibrations and weak equivalences are defined underlying in $\Top$, we see that $\mathcal{P} \mhyphen \mathbf{Alg}$ is a $\s\Set_{KQ}$-enriched model category with respect to this enrichment. By \autoref{example:InheritedCofibrationTestCategory}, any small full subcategory closed under finite colimits and tensors with finite simplicial sets inherits the structure of a cofibration test category.
\end{example}

\begin{example}
One can modify the previous example in a way that is similar to \autoref{example:ConvenientCofTestCatSpaces}. Namely, suppose that $\bfC \subset \mathcal{P} \mhyphen \mathbf{Alg}$ is a small full subcategory closed under finite colimits and tensors by finite simplicial sets, and suppose that $F|X|$ is contained in $\bfC$ for any finite simplicial set $X$, where $F \colon \Top \to \mathcal{P} \mhyphen \mathbf{Alg}$ is the left adjoint of the free-forgetful adjunction. Define the full subcategory of test objects $\bfT \subset \bfC$ to be the category of objects of the form $F|X|$ for $X$ a finite simplicial set, and define the (trivial) cofibrations to be the maps of the form $F|i| \colon F|X| \to F|Y|$ where $i$ is a (trivial) cofibration between finite simplicial sets in $\s\Set_{KQ}$. Then $(\bfC,\bfT)$ is a cofibration test category, hence we obtain a model structure on $\Ind(\bfC)$ by \autoref{theorem:ModelStructureIndC}. Since the inclusion $\bfC \hookrightarrow \mathcal{P}\mhyphen\mathbf{Alg}$ preserves finite colimits, it induces an adjunction $\Ind(\bfC) \rightleftarrows \mathcal{P}\mhyphen\mathbf{Alg}$. One can show that this adjunction is a Quillen equivalence.
\end{example}


\section{Example: a convenient model category of topological spaces}\label{sec:ConvenientCategorySpaces}

Throughout this section, let $\Top$ be a convenient category of spaces, such as $k$-spaces, compactly generated weak Hausdorff spaces or compactly generated Hausdorff spaces. Suppose that a small full subcategory $\bfC \subset \Top$ is given that is closed under finite colimits and tensors with finite simplicial sets, and that contains the space $|X|$ for any finite simplicial set $X$. As explained in \autoref{example:ConvenientCofTestCatSpaces}, if we define $\bfT$ to be the collection of spaces of the form $|X|$, where $X$ is any finite simplicial set, and if we define a map to be a (trivial) cofibration if and only if it is the geometric realization of a (trivial) cofibration in $\s\Set_{KQ}$ between finite simplicial sets, then $(\bfC,\bfT)$ is a cofibration test category. In this section, we will study this example in more detail.

We begin by characterizing the weak equivalences of $\Ind(\bfC)$. Note that the geometric realization functor $| \cdot | \colon \FinsSet \to \bfC$ extends uniquely to a filtered colimit preserving functor $| \cdot | \colon \s\Set \to \Ind(\bfC)$ that has a right adjoint $\SingTop$ defined by $(\SingTop C)_n = \Hom(|\Delta^n|, C)$ for any $C \in \Ind(\bfC)$.

\begin{lemma}\label{lemma:SingDetectsWEs}
Let $(\bfC,\bfT)$ be a cofibration test category as above. Then a map $C \to D$ in $\Ind(\bfC)$ is a weak equivalence if and only if $\Map(*,C) \to \Map(*,D)$ is a weak equivalence, where $*$ is the terminal object. In particular, $C \to D$ is a weak equivalence if and only if $\SingTop C \to \SingTop D$ is a weak equivalence in $\s\Set_{KQ}$.
\end{lemma}

\begin{proof}
	If $C \to D$ is a weak equivalence in $\Ind(\bfC)$, then $\Map(*,C) \to \Map(*,D)$ is a weak equivalence by definition. Conversely, suppose that $\Map(*,C) \to \Map(*,D)$ is a weak equivalence and let $X$ be a finite simplicial set. It follows by adjunction that $\Map(|X|,C) \to \Map(|X|,D)$ agrees with $\Map(*,C)^X \to \Map(*,D)^X$, hence this map is a weak equivalence.
	
	For the second statement, note that $\Map(*,E) \cong \Hom(* \otimes \Delta^\bullet, E) \cong \SingTop E$.
\end{proof}

The inclusion $\bfC \hookrightarrow \Top$ induces an adjunction $L : \Ind(\bfC) \rightleftarrows \Top : \wh{(\cdot)}_{\Ind}$ as explained in \autoref{ssec:ProAndIndCategories}, where $L$ is defined by $L(\{c_i\}) = \colim_i c_i$ for any $\{c_i\}$ in $\Ind(\bfC)$. Since geometric realization commutes with colimits, we see that the geometric realization functor $|\cdot| \colon \s\Set \to \Top$ factors as $\s\Set \xrightarrow{|\cdot|} \Ind(\bfC) \xrightarrow{L} \Top$.

\begin{proposition}\label{proposition:ConventientCategoryTopologicalEquivalences}
Let $(\bfC,\bfT)$ be a cofibration test category as above. The adjunctions $|\cdot| : \s\Set_{KQ} \rightleftarrows \Ind(\bfC) : \SingTop$ and $L : \Ind(\bfC) \rightleftarrows \Top : \wh{(\cdot)}_{\Ind}$ are Quillen equivalences.
\end{proposition}

\begin{proof}
It is clear from the definition of the (trivial) cofibrations in $(\bfC,\bfT)$ that $|\cdot| \colon \s\Set_{KQ} \to \Ind(\bfC)$ and $L \colon \Ind(\bfC) \to \Top$ send generating (trivial) cofibrations to (trivial) cofibrations. In particular, they are left Quillen functors. Since the composition of these adjunctions is the well-known Quillen equivalence $|\cdot| : \s\Set_{KQ} \rightleftarrows \Top : \SingTop$, it suffices to show by the 2 out of 3 property for Quillen equivalences that $|\cdot| : \s\Set_{KQ} \rightleftarrows \Ind(\bfC) : \SingTop$ is a Quillen equivalence. By \autoref{lemma:SingDetectsWEs}, a map $C \to D$ in $\Ind(\bfC)$ is a weak equivalence if and only if $\SingTop C \to \SingTop D$ is so. In particular, this adjunction is a Quillen equivalence if and only if the unit $X \to \SingTop |X|$ is a weak equivalence for any simplicial set $X$. If $X$ is a finite simplicial set, then $X \to \SingTop |X|$ agrees by definition with the unit of the adjunction $|\cdot| : \s\Set_{KQ} \rightleftarrows \Top : \SingTop$, which is always a weak equivalence. Since weak equivalences are stable under filtered colimits in $\s\Set_{KQ}$, it follows that the unit $X \to \SingTop |X|$ of $|\cdot| : \s\Set_{KQ} \rightleftarrows \Ind(\bfC) : \SingTop$ is a weak equivalence for any simplicial set $X$.
\end{proof}

One can show that the model category $\Ind(\bfC)$, with $(\bfC,\bfT)$ a cofibration test category of the type considered above, is very similar to $\Top$. We mention a few similarities. We first note that it is possible define homotopy groups for objects of $\Ind(\bfC)$, and that they detect weak equivalences. If $C$ is an object in $\Ind(\bfC)$, then by a basepoint of $C$ we mean a map $* \to C$.

\begin{definition}
The $n$-th homotopy group $\pi_n(C,c_0)$ of an object $C \in \Ind(\bfC)$ and a basepoint $c_0 \colon * \to C$ is defined as the set of pointed maps $|\Delta^n/\partial \Delta^n| \to C$ modulo pointed homotopy.
\end{definition}

It follows directly from this definition that $\pi_n(C,c_0) = \pi_n(\SingTop C,c_0)$ for any $C \in \Ind(\bfC)$ and $c_0 \in C$. We conclude the following:

\begin{proposition}\label{Proposition:Convenient1}
A map $f \colon C \to D$ in $\Ind(\bfC)$ is a weak equivalence if and only if $\pi_n(C,c_0) \to \pi_n(D,f(c_0))$ is a bijection for any $c_0 \in C$ and $n \geq 0$. Moreover, the homotopy groups for objects in $\Ind(\bfC)$ commute with filtered colimits.
\end{proposition}

\begin{proof}
The first statement follows since $f$ is a weak equivalence if and only if $\SingTop f$ is so. The second part follows since both the functor $\SingTop$ and the homotopy groups of simplicial sets commute with filtered colimits.
\end{proof}

It can also be shown that one can take the same generating (trivial) cofibrations in $\Ind(\bfC)$ as in the usual Quillen model structure on $\Top$. Define
\[I = \{ \partial D^n \hookrightarrow D^n \mid n \geq 0 \} \]
and
\[J = \{ D^n \times \{0\} \hookrightarrow D^n \times [0,1] \mid n \geq 0\} .\]

\begin{proposition}\label{Proposition:Convenient2}
The sets $I$ and $J$ are sets of generating cofibrations and generating trivial cofibrations for $\Ind(\bfC)$, respectively.
\end{proposition}

\begin{proof}
We need to show that the geometric realization of any cofibration (or trivial cofibration) between finite simplicial sets lies in the saturation of $I$ (or $J$, respectively). This follows from the fact that $|\cdot| \colon \s\Set \to \Ind(\bfC)$ preserves colimits and that all maps of the form $|\partial \Delta^n| \to |\Delta^n|$ (or $|\Lambda^n_k| \to |\Delta^n|$) are isomorphic to a map in $I$ (or $J$, respectively).
\end{proof}

For any two objects $C = \{c_i\}_i$ and $D = \{d_j\}_j$ in $\Ind(\bfC)$, one can compute their product levelwise by $C \times D = \{c_i \times d_j\}_{(i,j) \in I \times J}$. Since the finite colimits of $\bfC$ are computed in $\Top$, and finite colimits in $\Ind(\bfC)$ can be computed levelwise, we see that the functor $- \times D \colon \bfC \to \Ind(\bfC)$ preserves finite colimits for any $D$ in $\Ind(\bfC)$. As explained in \autoref{ssec:ProAndIndCategories}, it follows from this that the product functor $- \times D \colon \Ind(\bfC) \to \Ind(\bfC)$ has a right adjoint. In particular, $\Ind(\bfC)$ is cartesian closed. This cartesian closed structure interacts well with the model structure defined above.

\begin{proposition}\label{Proposition:Convenient3}
$\Ind(\bfC)$ is a cartesian closed model category.
\end{proposition}

\begin{proof}
It suffices to show that for any pair of generating cofibrations $C \cofarrow D$, $C' \cofarrow D'$, the pushout-product
\[C \times D' \cup_{C \times C'} D \times C' \to D \times D' \]
is a cofibration that is trivial if either $C \cofarrow D$ or $C' \cofarrow D'$ is so. This is clearly true.
\end{proof}

One can furthermore show that the full subcategory of $\Top$ on the CW-complexes embeds fully faithfully into $\Ind(\bfC)$. Note that any finite CW-complex $X$ is (homeomorphic to) an object in $\bfC$.

\begin{proposition}\label{Proposition:Convenient4}
There is a fully faithful functor from the category of CW-complexes into $\Ind(\bfC)$ that preserves and detects weak equivalences.
\end{proposition}

\begin{proof}
If $X$ is a CW-complex, then one can always choose a CW-decomposition. The finite CW-subcomplexes in this decomposition together with their inclusions form a directed diagram $\{X_i\}$ for which $\colim_i X_i \cong X$. Suppose that we have chosen a CW-decomposition for any CW-complex $X$ and denote the associated directed diagram of finite CW-subcomplexes by $\{X_{i_X}\}$. Since a map from a compact space into a CW-complex (with a given CW-decomposition) always lands in a finite CW-subcomplex, we see that the canonical map 
\[\lim_{i_X} \colim_{i_Y} \Hom(X_{i_X}, Y_{i_Y}) \to \Hom(X,Y)\]
is an isomorphism for any pair of CW-complexes. By definition of the morphisms in $\Ind(\bfC)$, this implies that the functor that sends a CW-complex $X$ to the ind-object $\{X_{i_X}\}$ in $\Ind(\bfC)$ is well-defined and fully faithful. Preservation and detection of weak equivalences follows directly form the fact that $\SingTop$ detects weak equivalences and that $\colim_{i_X} \SingTop(X_{i_X}) \cong \SingTop(X)$ for any CW-complex $X$.
\end{proof}

We end this section by discussing a specific example of such a full subcategory $\bfC$, namely the category $\mathbf{CM}$ of compact metrizable spaces. Under Gelfand-Naimark duality, this category corresponds to the category of separable commutative unital $C^*$-algebras. If we let $\Top$ be the category of compactly generated Hausdorff spaces, then $\mathbf{CM}$ as a full subcategory is closed under all finite colimits and tensors by finite simplicial sets. In particular, by the above we obtain a model structure on $\Ind(\mathbf{CM})$ that is equivalent to the Quillen model structure on $\Top$. In \cite{Barnea2017Convenient}, Barnea also proposes a model structure on $\Ind(\mathbf{CM})$. However, this model structure does not agree with the one constructed above, so we will briefly describe his model structure and the difference with ours. Let us denote our model structure by $\Ind(\mathbf{CM})_Q$.

Barnea shows in \cite{Barnea2017Convenient} that $\mathbf{CM}$ is a ``special weak cofibration category'', and hence that there exists an induced model structure on $\Ind(\mathbf{CM})$, which we will denote by $\Ind(\mathbf{CM})_B$. This model structure is cofibrantly generated and one can take the set of Hurewicz cofibrations in $\mathbf{CM}$ as a set of generating cofibrations, while one can take the Hurewicz cofibrations that are also homotopy equivalences as a set of generating trivial cofibrations. If we define $\bfT = \mathbf{CM}$ and if we define a map in $\bfT$ to be a (trivial) cofibration if it is in the set of generating (trivial) cofibrations just mentioned, then $(\mathbf{CM},\bfT)$ is a cofibration test category and the completed model structure on $\Ind(\mathbf{CM})$ coincides with the one that Barnea constructed. Since Barnea's model structure $\Ind(\mathbf{CM})_B$ has strictly more generating (trivial) cofibrations than our model structure $\Ind(\mathbf{CM})_Q$, we see that the identity functor is a left Quillen functor $\Ind(\mathbf{CM})_Q \to \Ind(\mathbf{CM})_B$. To see that the model structures do not coincide, we will show that $\Ind(\mathbf{CM})_Q$ has strictly more weak equivalences than $\Ind(\mathbf{CM})_B$. Let $C$ be any metrizable infinite Stone space, such as a Cantor space. Then, for $\SingTop$ and $|\cdot|$ as defined just above \autoref{lemma:SingDetectsWEs}, the counit $|\SingTop C| \to C$ is a weak equivalence in $\Ind(\mathbf{CM})_Q$. However, this map is not a weak equivalence in $\Ind(\mathbf{CM})_B$, since $\Map(C,|\SingTop C|) \to \Map(C,C)$ is not a weak equivalence of simplicial sets: Since these mapping spaces are discrete, this would imply that the map is an isomorphism. However, it is not surjective since there is no map $C \to |\SingTop C|$ that gets mapped to $\id_C$. The model structure $\Ind(\mathbf{CM})_Q$ defined here is similar to the Quillen model structure on $\Top$, while Barnea's model structure $\Ind(\mathbf{CM})_B$ bears some similarity to the Strøm model structure on $\Top$.


\section{The dual model structure on \texorpdfstring{$\Pro(\bfC)$}{Pro(C)}}\label{sec:ProCompletedModelStructure}
A model structure on $\bfE$ also gives rise to a model structure on $\bfE^{op}$, where the fibrations (cofibrations) of $\bfE^{op}$ are the cofibrations (fibrations, respectively) of $\bfE$. In particular, $\bfE$ is cofibrantly generated if and only if $\bfE^{op}$ is fibrantly generated. Since $\Pro(\bfC) \simeq \Ind(\bfC^{op})^{op}$, this implies that if $\bfC$ is the dual of a cofibration test category, then $\Pro(\bfC)$ admits a fibrantly generated simplicial model structure. We explicitly dualize the main definition and result of \autoref{sec:IndCompletedModelStructure} in this section, and then discuss a few examples of such fibrantly generated simplicial model structures on pro-categories. Again, we work with $\s\Set$ endowed with either the Joyal or the Kan-Quillen model structure.

We say that a simplicial category $\bfC$ is \emph{finitely cotensored} if $\bfC^{op}$ is finitely tensored in the sense of \autoref{definition:FinitelyTensored}. Explicitly, this means that $\bfC$ admits finite limits and cotensors by finite simplicial sets, and that these commute with each other. As explained in \autoref{ssec:ProAndIndCategories}, if $\bfC$ is a small simplicial category that is finitely cotensored, then the simplicial category $\Pro(\bfC)$ is tensored, cotensored, complete and cocomplete.

\begin{definition}\label{definition:FibTestCat}
A \emph{fibration test category} $(\bfC, \bfT)$ consists of a small finitely cotensored simplicial category $\bfC$, a full subcategory $\bfT \subset \bfC$ of \emph{test objects} and two classes of maps in $\bfT$ called \emph{fibrations}, denoted $\fibarrow$, and \emph{trivial fibrations}, denoted $\trivfibarrow$, both containing all isomorphisms, that satisfy the following properties:
\begin{enumerate}[(1)]
    \item\label{FibTestCat:TestIsFibrant} The terminal object $*$ is a test object, and for every test object $t \in \bfT$, the map $t \to *$ is a fibration.
    \item\label{FibTestCat:TestObjectsClosedUnderCotensors} For every fibration between test objects $s \fibarrow t$ and cofibration between finite simplicial set $U \cofarrow V$, the pullback-power map $s^V \to s^U \times_{t^U} t^V$ is a fibration between test objects which is trivial if either $s \fibarrow t$ or $U \cofarrow V$ is so.
    \item\label{FibTestCat:TestingTrivialFibration} A morphism $c \to d$ in $\bfT$ is a trivial fibration if and only if it is a fibration and $\Map(d,t) \to \Map(c,t)$ is a weak equivalence of simplicial sets for every $t \in \bfT$.
    \item\label{FibTestCat:FibrantMappingSpaces} Any object $c \in \bfC$ has the left lifting property with respect to trivial fibrations.
\end{enumerate}
\end{definition}

For a fibration test category $\bfC$, we write $\fib(\bfC)$ for the set of fibrations and $\we(\bfC)$ for the set of maps $c \to d$ that induce a weak equivalence $\Map(d,t) \to \Map(c,t)$ for every $t \in \bfT$. By property \eqref{FibTestCat:TestingTrivialFibration}, the set of trivial fibrations is $\fib(\bfC) \cap \we(\bfC)$. Note that the definition of a fibration test category is formally dual to that of a cofibration test category. More precisely, $(\bfC,\bfT)$ is a fibration test category if and only if $(\bfC^{op}, \bfT^{op})$ is a cofibration test category in the sense of \autoref{definition:CofTestCat}, where the (trivial) cofibrations of $(\bfC^{op}, \bfT^{op})$ are defined as the (trivial) fibrations of $(\bfC,\bfT)$.

We let $P \subset \Ar(\Pro(\bfC))$ denote the image of the set $\fib(\bfC)$ along the inclusion $\bfC \hookrightarrow \Pro(\bfC)$, and $Q \subset \Ar(\Pro(\bfC))$ the image of the set of trivial fibrations. The sets $P$ and $Q$ are the generating (trivial) fibrations of the completed model structure on $\Pro(\bfC)$. The following theorem is formally dual to \autoref{theorem:ModelStructureIndC}.

\begin{theorem}\label{theorem:ModelStructureProC}
Let $(\bfC,\bfT)$ be a fibration test category. Then $\Pro(\bfC)$ carries a fibrantly generated (hence cocombinatorial) simplicial model structure, the \emph{completed model structure}, where a map $C \to D$ is a weak equivalence if and only if $\Map(D,t) \to \Map(C,t)$ is a weak equivalence for every $t \in \bfT$. A set of generating fibrations (generating trivial fibrations) is given by $P$ ($Q$, respectively). Every object is cofibrant in this model structure.
\end{theorem}

\begin{example}\label{example:InheritedFibrationTestCategory}
Dualizing \autoref{example:InheritedCofibrationTestCategory}, we see that if $\bfE$ is a simplicial model category in which every object is cofibrant, then any small full subcategory $\bfC \subset \bfE$ closed under finite limits and finite cotensors admits the structure of a fibration test category. Namely, defining $\bfT$ to be the full subcategory of fibrant objects of $\bfC$, and defining the (trivial) fibrations to be those of $\bfE$ between objects in $\bfT$, then $(\bfC,\bfT)$ is a fibration test category. As in \autoref{example:InheritedCofibrationTestCategory}, will say that $\bfC$ \emph{inherits} the structure of a fibration test category from $\bfE$.
\end{example}

\begin{remark}
Observe that for the fibration test category $(\bfC,\bfT)$ from the previous example, the completed model structure on $\Pro(\bfC)$ is a special case of \autoref{thm:ApproximateTheorem}, namely the case where $\bfT$ is the collection of all fibrant objects in $\bfC$. By (the dual of) \autoref{example:RestrictingToLessTestObjects}, it follows that we can take $\bfT$ to be any collection of fibrant objects in $\bfC$ that is closed under ``finite pullback-powers''. The general case, where we let $\bfT$ be any collection of fibrant objects in $\bfC$, is discussed in \autoref{sec:BousfieldLocalizations}.
\end{remark}

\begin{example}\label{example:KanQuillenLeanFibrationTestCategory}
Recall that we call a simplicial set \emph{lean} if it is degreewise finite and coskeletal. The full subcategory of $\s\Set$ spanned by all lean simplicial sets $\LsSet$ is closed under finite limits and finite cotensors. By \autoref{example:InheritedFibrationTestCategory}, it inherits the structure of a fibration test category from $\s\Set_{KQ}$, which we will denote by $\LsSet_{KQ}$. By \autoref{theorem:ModelStructureProC} we obtain a model structure on $\Pro(\LsSet)$. Since this category is equivalent to the category of simplicial profinite sets $\s\wh\Set$ by (the dual of) \autoref{thm:IndValuedPresheavesIsIndSkeletalPresheaves}, we in particular obtain a simplicial model structure on $\s\wh\Set$. This model structure coincides with Quick's model structure for profinite spaces \cite{Quick2008Profinite}, as explained in \autoref{corollary:QuicksModelStructureAgrees} below. We denote it by $\s\wh\Set_Q$.
\end{example}

\begin{example}\label{example:MorelsModelStructures}
Consider the full simplicial subcategory $\bfT_p$ of $\s\Set$ whose objects are those lean Kan complexes that have finite $p$-groups as homotopy groups. One can show that $\bfT_p$ is closed under ``finite pullback-powers'', so by the previous example and the dual of \autoref{example:RestrictingToLessTestObjects}, we obtain a fibration test category $\LsSet_p = (\LsSet,\bfT_p)$ in which the (trivial) fibrations are the (trivial) Kan fibrations between objects of $\bfT_p$. It is proved in \autoref{corollary:MorelsModelStructureAgrees} that the completed model structure on $\Pro(\LsSet_p)$ agrees with Morel's model structure for pro-$p$ spaces \cite{Morel1996EnsemblesProfinis}.
\end{example}

\begin{example}\label{example:JoyalFibrationTestCategory}
The category of lean simplicial sets also inherits the structure of a fibration test category from the Joyal model structure $\s\Set_J$, which we will denote by $\LsSet_J$. The corresponding model structure on $\s\wh\Set$ obtained from \autoref{theorem:ModelStructureProC} will be called the \emph{profinite Joyal model structure} and its fibrant objects will be called \emph{profinite quasi-categories}. We will come back to this model category in \autoref{sec:CompleteSegalProfiniteVSProfiniteQuasiCats}, and we will describe its underlying $\infty$-category in \autoref{remark:UnderlyingInftyCategoryProfiniteInftyCats}.
\end{example}

\begin{example}\label{example:PStratifiedFibrationTestCategory}
In \cite{Haine2018HomotopyV5}, Haine defines the Joyal-Kan model structure on $\s\Set_{/P}$, where $P$ is (the nerve of) a poset. This model category describes the homotopy theory of $P$-stratified spaces. Since it is a left Bousfield localization of the Joyal model structure on $\s\Set_{/P}$, any object is cofibrant and it is a $\s\Set_J$-enriched model category. Actually, this model structure can be shown to be $\s\Set_{KQ}$-enriched \cite[\S\S 2.4-5]{Haine2018HomotopyV5}. In particular, any small full subcategory $\bfC$ closed under finite limits and cotensors by finite simplicial sets inherits the structure of a fibration test category. If $P$ is a finite poset and $\bfC = \LsSet_{/P}$ is the full subcategory of lean simplicial sets over (the nerve of) $P$, then one can show that $\Pro(\LsSet_{/P}) \cong \s\wh\Set_{/P}$. In particular, by \autoref{theorem:ModelStructureProC}, we obtain a model structure on $\s\wh\Set_{/P}$ that is $\s\Set_{KQ}$-enriched. It is shown in \autoref{example:UnderlyingInftyCatProfiniteStratified} that the underlying $\infty$-category of this model category is the $\infty$-category of profinite $P$-stratified spaces defined in \cite{BarwickGlasmanHaine2018ExodromyV7}.
\end{example}

\begin{example}\label{example:ProfiniteGroupoids}
We call a groupoid finite if it has finitely many arrows (including the identity arrows). The category of finite groupoids $\Fin\Grpd$ inherits the structure of a fibration test category from the canonical model structure on $\Grpd$ \cite[\S 5]{Anderson1978FibrationsRealization}. (Note that $\Grpd$ can be viewed as a $\s\Set_{KQ}$-enriched model structure by defining $\Map(A,B) = N(\Fun(A,B))$ for any $A,B \in \Grpd$.) The completed model structure on the category of profinite groupoids $\wh \Grpd = \Pro(\Fin\Grpd)$ obtained from \autoref{theorem:ModelStructureProC} coincides with the model structure for profinite groupoids defined by Horel in \cite[\S 4]{Horel2017ProfiniteOperads}. To see this, note that Horel shows in \cite[\S 4]{Horel2017ProfiniteOperads} that the Barnea-Schlank model structure on $\wh\Grpd$ exists and coincides with his model structure. By \autoref{remark:BarneaSchlankVsFibrationTestCategory} below, the Barnea-Schlank model structure on $\wh\Grpd$ must coincide with our model structure. In particular, Horel's model structure agrees with the one that we construct in this example.
\end{example}

\begin{example}
Similarly, we call a category finite if it has finitely many arrows. The category of all categories admits the canonical model structure, defined for example in \cite{Rezk2000ModelCategoryOfCategories}. Since this model structure is $\s\Set_J$-enriched, the category of finite categories $\Fin\Cat$ inherits the structure of a fibration test category. By \autoref{theorem:ModelStructureProC}, we obtain a $\s\Set_J$-enriched model structure on $\wh \Cat = \Pro(\Fin\Cat)$ which we will call the \emph{model structure for profinite categories}.
\end{example}

\begin{example}\label{example:ReedyFibrationTestCategory}
Let $\bisSet$ be endowed with the Reedy model structure with respect to the Kan-Quillen model structure on $\s\Set$. Recall that the category of bisimplicial profinite sets $\biswhSet$ is equivalent to $\Pro(\LsSet^{(2)})$, where $\LsSet^{(2)}$ denotes the category of doubly lean bisimplicial sets defined at the end of \autoref{ssec:ProAndIndCategories}. Since any object in $\bisSet$ is cofibrant, $\LsSet^{(2)}$ inherits the structure of a fibration test category from the Reedy model structure on $\bisSet$. By applying \autoref{theorem:ModelStructureProC}, we obtain a model structure on $\Pro(\LsSet^{(2)}) \simeq \biswhSet$. This model structure coincides with the Reedy model structure on $\biswhSet$ with respect to the Quick model structure on $\s\wh\Set$, as will be shown in \autoref{proposition:ReedyModelStructuresCoincide}.
\end{example}

\begin{remark}\label{remark:BarneaSchlankVsFibrationTestCategory}
As discussed in the introduction, there are similarities between our construction of a model structure on $\Pro(\bfC)$ and the construction of Barnea-Schlank in \cite{BarneaSchlank2016ProSimplicialSheaves}. Suppose $\bfC$ is a fibration test category in the sense of \autoref{definition:FibTestCat}. Then $\bfC$ comes with a set $\fib(\bfC)$ of fibrations and a set of $\we(\bfC)$ of weak equivalences. It is very unlikely that the triple $(\bfC,\fib(\bfC),\we(\bfC))$ is a ``weak fibration category'' in the sense of Definition 1.2 of \cite{BarneaSchlank2016ProSimplicialSheaves}. Namely, that definition asks that $\fib(\bfC)$ contains all isomorphisms of $\bfC$, that it is closed under composition, and that a pushout of a map in $\fib(\bfC)$ is again in $\fib(\bfC)$. However, if we define $\fib'(\bfC)$ to be the smallest set that contains $\fib(\bfC)$ and that satisfies these properties, then $(\bfC,\fib'(\bfC),\we(\bfC))$ might be a weak fibration category. If this is the case, then the ``induced'' model structure on $\Pro(\bfC)$, in the sense of Theorem 1.8 of \cite{BarneaSchlank2016ProSimplicialSheaves}, could exist. The cofibrations of this model structure are defined as the maps that have the left lifting property with respect to $\fib'(\bfC) \cap \we(\bfC)$, while the trivial cofibrations are the maps that have the left lifting property with respect to $\fib'(\bfC)$. Since the maps in $\fib'(\bfC)$ are clearly fibrations in our construction of the ``completed model structure'' on $\Pro(\bfC)$ (see \autoref{theorem:ModelStructureProC}), we conclude that the (trivial) cofibrations for both model structures must agree. In particular, if both our model structure and the Barnea-Schlank model structure of \cite{BarneaSchlank2016ProSimplicialSheaves} exist on $\Pro(\bfC)$, then they must coincide. An example where this happens is when $\bfC = \Fin\Grpd$. (See \autoref{example:ProfiniteGroupoids} above.)
\end{remark}


\section{Comparison to some known model structures}\label{sec:Comparison}
As stated in \autoref{theorem:ModelStructureProC}, for any fibration test category $(\bfC,\bfT)$, all objects in the completed model structure on $\Pro(\bfC)$ are cofibrant. We will now show that, in the case that $\bfC$ is the category $\LsSet$ of lean simplicial sets, this statement can often be strengthened to say that the cofibrations are exactly the monomorphisms. We show how this can be used to prove that the model structures on $\s\wh\Set$ obtained in \autoref{example:KanQuillenLeanFibrationTestCategory} and \autoref{example:MorelsModelStructures} agree with Quick's model structure and Morel's model structure, respectively. It will also follow that the cofibrations in the profinite Joyal model structure from \autoref{example:JoyalFibrationTestCategory} are exactly the monomorphisms. We conclude this section by showing that the model structure on $\biswhSet$ from \autoref{example:ReedyFibrationTestCategory} agrees with the Reedy model structure on $\biswhSet$ with respect to Quick's model structure on $\s\wh\Set$. The main result about cofibrations in $\s\wh\Set$ is the following:

\begin{proposition}\label{appendix:proposition:CofibrationsAreMonosSimplicial}
Let $\LsSet$ be the category of lean simplicial sets endowed with the structure of a fibration test category. Suppose that for any contractible lean Kan complex $K$, the map $K \to *$ is a trivial fibration in $\LsSet$, and further that any trivial fibration $L \trivfibarrow K$ in $\LsSet$ is a trivial Kan fibration. Then the cofibrations in the completed model structure on $\Pro(\LsSet) \simeq \s\wh\Set$ are the monomorphisms.
\end{proposition}

This proposition clearly applies to the fibration test categories $\LsSet_{KQ}$, $\LsSet_p$ and $\LsSet_J$ of \Cref{example:KanQuillenLeanFibrationTestCategory,example:MorelsModelStructures,example:JoyalFibrationTestCategory}. The following lemmas will be used to prove this result. Recall that the category of profinite sets $\wh \Set$ is equivalent to the category of Stone spaces $\Stone$.

\begin{lemma}\label{appendix:lemma:LevelRepresentationByMonos}
A map of profinite sets (or simplicial profinite sets) $S \to T$ is a monomorphism if and only if it is (isomorphic to) the limit of a cofiltered diagram $\{S_i \cofarrow T_i\}_{i \in I}$ consisting of monomorphisms between finite sets (or degreewise finite simplicial sets, respectively). 
\end{lemma}

\begin{proof}
Note that in the category of Stone spaces, the monomorphisms are precisely the injective continuous maps. Since a cofiltered limit of injective maps is again injective, we see that if $S \to T$ is an inverse limit of monomorphisms $S_i \cofarrow T_i$, then $S \to T$ is itself a monomorphism.

Conversely, suppose that $S \to T$ is a monomorphism of profinite sets (or simplicial profinite sets). Write $T = \lim_i T_i$ as a cofiltered limit of finite sets (or lean simplicial sets, respectively), and, for every $i$, write $S'_i$ for the image of the composition $S \to T \to T_i$. Then $\{S'_i\}_{i \in I}$ is a cofiltered diagram since the structure maps $T_i \to T_j$ restrict to maps $S'_i \to S'_j$ for any $i \to j$ in $I$. Since $\{S_i' \to T_i\}_{i \in I}$ is levelwise a monomorphism, the proof is complete if we can show that $S \to \lim_i S_i'$ is an isomorphism. Since isomorphisms of Stone spaces are detected on the underlying sets, it suffices to show that this map is both injective and surjective. It is injective since the composition $S \to \lim_i S'_i \to T$ is, while it is surjective by \cite[Corollary 1.1.6]{RibesZalesskii2010ProfiniteGroups}
\end{proof}

We will denote the two-element set $\{0,1\}$ by $\bftwo$.

\begin{lemma}\label{appendix:lemma:RecognizingMonoProfiniteSets}
A map of (profinite) sets $S \to T$ is a monomorphism if and only if it has the left lifting property with respect to $\bftwo \to *$.
\end{lemma}

\begin{proof}
We leave the case where $S \to T$ is a map of sets to the reader. For the ``if'' direction in the profinite case, suppose that $f \colon S \to T$ has the left lifting property with respect to $\bftwo \to *$, but is not a monomorphism. Regarding $S$ and $T$ as Stone spaces, there must exist distinct $s,s' \in S$ such that $f(s) = f(s')$. Choose some clopen $U \subset S$ such that $s \in U$ and $s' \not \in U$. Then the indicator function $\mathbb{1}_U \colon S \to \bftwo$ is continuous but does not extend to a map $T \to \bftwo$. We conclude that $S \to T$ must be a monomorphism.

For the converse, note that by \autoref{appendix:lemma:LevelRepresentationByMonos} we may assume without loss of generality that $S \to T$ can be represented by levelwise monomorphisms $\{S_i \to T_i\}$. Since $\bftwo$ is cocompact in $\wh\Set$, any map $S \to \bftwo$ factors through $S_i$ for some $i$. Since $S_i \to T_i$ is a monomorphism of sets, the result follows.
\end{proof}

Consider the diagram
\[\begin{tikzcd}
* \ar[r,"\bftwo"] \ar[d,"{[n]}"'] & \Fin\Set \\
\Delta^{op} \ar[ur,dashed,"R_n \bftwo"'] & 
\end{tikzcd}\]
where $[n]$ denotes the inclusion of the terminal category $*$ into $\Delta^{op}$ at $[n]$, and where $\bftwo$ denotes the inclusion of $*$ into $\Fin\Set$ at the two-element set $\bftwo$. Since $\Fin\Set$ has all finite limits, the right Kan extension $R_n \bftwo$ exists. Since the inclusion $* \hookrightarrow \Delta^{op}$ factors through $\Delta_{\leq n}^{op}$, the simplicial set $R_n \bftwo$ is $n$-coskeletal. In particular, it is a lean simplicial set.

\begin{lemma}\label{appendix:lemma:RecognizingMonoSimplicialProfiniteSets}
A map of simplicial (profinite) sets is a monomorphism if and only if it has the left lifting property with respect to $R_n \bftwo \to *$ for every $n \in \bbN$.
\end{lemma}

\begin{proof}
Since the inclusions of $\Fin\Set$ into $\Set$ and $\wh\Set$ both preserve limits, we see that the lean simplicial set $R_n \bftwo$ constructed above is also the right Kan extension of $* \xrightarrow{[n]} \Delta^{op}$ along $* \xrightarrow{\bftwo} \Set$ and along $* \xrightarrow{\bftwo} \wh\Set$. In particular, a map of simplicial (profinite) sets $X \to Y$ has the left lifting property with respect to $R_n \bftwo \to *$ if and only if $X_n \to Y_n$ has the left lifting property with respect to $\bftwo \to *$, hence the result follows from \autoref{appendix:lemma:RecognizingMonoProfiniteSets}.
\end{proof}

\begin{proof}[Proof of \autoref{appendix:proposition:CofibrationsAreMonosSimplicial}]
We first show that any cofibration in the model category $\Pro(\LsSet)$ is a monomorphism. Since $R_n\bftwo \to *$ has the right lifting property with respect to all monomorphisms in $\s\Set$, we see that it is a trivial Kan fibration, hence by assumption a trivial fibration in the fibration test category $\LsSet$ and a generating trivial fibration in $\Pro(\LsSet)$. By \autoref{appendix:lemma:RecognizingMonoSimplicialProfiniteSets}, any cofibration in $\Pro(\LsSet) \simeq \s\wh\Set$ is a monomorphism.

For the converse, suppose $X \to Y$ is a monomorphism in $\Pro(\LsSet)$. By \autoref{appendix:lemma:LevelRepresentationByMonos}, we may assume that $X \to Y$ is a cofiltered limit of monomorphisms between degreewise finite simplicial sets $\{X_i \to Y_i\}_{i \in I}$. We see that for every $i$, the map $X_i \to Y_i$ has the left lifting property with respect to the generating trivial fibrations of $\Pro(\LsSet)$ since these are trivial Kan fibrations between lean simplicial sets. Since any generating trivial fibration is a map between cocompact objects, it follows that $X \to Y$ also has the left lifting property with respect to the generating trivial fibrations.
\end{proof}

\autoref{appendix:proposition:CofibrationsAreMonosSimplicial} shows that, for $\LsSet_{KQ}$ and $\LsSet_p$ the fibration test categories of \Cref{example:KanQuillenLeanFibrationTestCategory,example:MorelsModelStructures}, the cofibrations of the model categories $\Pro(\LsSet_{KQ})$ and $\Pro(\LsSet_p)$ are the monomorphisms. This means that the cofibrations coincide with those of Quick's model structure \cite{Quick2008Profinite} and Morel's model structure \cite{Morel1996EnsemblesProfinis}, respectively. The same is true for the weak equivalences. This follows from the results in \S 7 of \cite{BarneaHarpazHorel2017} (most notably Lemma 7.4.7 and Lemma 7.4.10), using that Quick's and Morel's model structures on $\s\wh\Set$ are simplicial model structures.

\begin{proposition}[\cite{BarneaHarpazHorel2017}] \label{prop:BHH17}
A map $X \to Y$ of simplicial profinite sets is a weak equivalence in Quick's model structure if and only if $\Map(Y,K) \to \Map(X,K)$ is a weak equivalence for any lean Kan complex $K$. It is a weak equivalence in Morel's model structure if and only if $\Map(Y,K) \to \Map(X,K)$ is a weak equivalence for any lean Kan complex $K$ whose homotopy groups are finite $p$-groups.
\end{proposition}

From this proposition and the definition of the completed model structure (see \autoref{theorem:ModelStructureProC}), we see that the weak equivalences of $\Pro(\LsSet_{KQ})$ (or $\Pro(\LsSet_p)$) agree with the weak equivalences in Quick's model structure (or Morel's model structure, respectively) on $\s\wh\Set$.

\begin{corollary}\label{corollary:QuicksModelStructureAgrees}
The completed model structure on $\Pro(\LsSet_{KQ})$ coincides with Quick's model structure.
\end{corollary}

\begin{corollary}\label{corollary:MorelsModelStructureAgrees}
For any prime number $p$, the completed model structure on $\Pro(\LsSet_p)$ coincides with Morel's model structure.
\end{corollary}

The proof of \autoref{appendix:proposition:CofibrationsAreMonosSimplicial} admits an analogue for bisimplicial sets (in fact, for the category of presheaves on $K$ for any small category $K$ that can be written as a union of finite full subcategories), which we leave as an exercise to the reader.

\begin{proposition}\label{appendix:proposition:CofibrationsAreMonosBisimplicial}
Let $\bisSet$ be endowed with a simplicial model structure in which the cofibrations are the monomorphisms, and let $\LsSet^{(2)}$ be the full subcategory of doubly lean bisimplicial sets, which inherits the structure of a fibration test category in the sense of \autoref{example:InheritedFibrationTestCategory}. Then the cofibrations in $\Pro(\LsSet^{(2)}) \simeq \biswhSet$ are the monomorphisms.
\end{proposition}

Note that this proposition implies that the cofibrations in the model structure on $\biswhSet$ from \autoref{example:ReedyFibrationTestCategory} are exactly the monomorphisms. We will show that, in fact, this model structure coincides with the Reedy model structure on $\biswhSet$ with respect to $\s\wh\Set_Q$. We do this by inspecting the generating (trivial) fibrations of the Reedy model structure. For the following proof, note that Quick's model structure coincides with the completed model structure on $\Pro(\LsSet_{KQ})$ by \autoref{corollary:QuicksModelStructureAgrees}, hence that any (trivial) Kan fibration between lean Kan complexes is a (trivial) fibration in Quick's model structure.

\begin{proposition}\label{proposition:ReedyModelStructuresCoincide}
The model structure on $\Pro(\LsSet^{(2)})$ of \autoref{example:ReedyFibrationTestCategory} coincides with the Reedy model structure on $\biswhSet$ (with respect to $\s\wh\Set_Q$).
\end{proposition}

\begin{proof}
Note that if $L \to K$ is a (trivial) Reedy fibration between Reedy fibrant doubly lean bisimplicial sets, then $L_n$ and $M_n L \times_{M_n K} K_n$ are lean Kan complexes for every $n$. In particular, the map $L_n \to M_n L \times_{M_n K} K_n$ is a (trivial) fibration between lean Kan complexes for every $n$. This shows that any generating (trivial) fibration in $\Pro(\LsSet^{(2)})$ is a (trivial) fibration in the Reedy model structure on $\biswhSet$.

For the converse, note that the Reedy model structure on $\biswhSet$ is fibrantly generated. Its generating (trivial) fibrations are maps of the form
\begin{equation}\label{map:ReedyGeneratingFibration}
    \bfG_n L \to \partial \bfG_n L \times_{\partial \bfG_n K} \bfG_n K
\end{equation}
for any $n \geq 0$, where $\bfG_n$ is the right adjoint to the functor $X \mapsto X_n$, where $\partial \bfG_n$ is the right adjoint to latching object functor $X \mapsto L_n X$, and where $L \to K$ is a generating (trivial) fibration in $\s\wh\Set$. It can be shown using the right adjointness of $\bfG_n$ and $\partial \bfG_n$ that these functors restrict to functors $\LsSet \to \LsSet^{(2)}$. One can furthermore deduce from the adjointness that if $L$ and $K$ are fibrant in $\s\Set$, then both the domain and codomain of the map \eqref{map:ReedyGeneratingFibration} are Reedy fibrant in $\biswhSet$ and hence in $\bisSet$. This shows that any map of the form \eqref{map:ReedyGeneratingFibration}, with $L \to K$ a (trivial) fibration in $\LsSet$, is a (trivial) fibration in $\LsSet^{(2)}$. In particular, any generating (trivial) fibration in the Reedy model structure on $\biswhSet$ is a (trivial) fibration in $\Pro(\LsSet^{(2)})$. We conclude that both model structures coincide.
\end{proof}


\section{Quillen pairs}\label{sec:QuillenPairs}

As explained in \autoref{ssec:ProAndIndCategories}, there is an easy criterion for constructing adjunctions between ind-categories: if $\bfC$ is a small category that admits finite colimits and if $\bfE$ is any cocomplete category, then a functor $F \colon \Ind(\bfC) \to \bfE$ has a right adjoint if and only it preserves all colimits. Furthermore, these functors correspond to functors $\bfC \to \bfE$ that preserve all finite colimits. There is a dual criterion for pro-categories. In the simplicial case, this can be strengthened as in the following lemma.

If $\bfE$ is a tensored cocomplete simplicial category, then we say that colimits and tensors commute in $\bfE$ if the analogue of item \eqref{item3:FinitelyTensored} of \autoref{definition:FinitelyTensored} holds for all diagrams in $\bfE$ and all simplicial sets.

\begin{lemma}\label{lem:TensorPreservingHasEnrichedRightAdjoint}
Let $\bfC$ be a small finitely tensored simplicial category and let $\bfE$ be a tensored cocomplete simplicial category in which colimits and tensors commute. Then any simplicial functor $F \colon \bfC \to \bfE$ that preserves finite colimits and tensors with finite simplicial sets extends to a functor $\wt F \colon \Ind(\bfC) \to \bfE$ that admits a right adjoint. Moreover, this adjunction is an enriched adjunction.
\end{lemma}

\begin{proof}
The simplicial functor $\wt F \colon \Ind(\bfC) \to \bfE$ is defined on objects by $\wt F(\{c_i\}) = \colim_i F(c_i)$ and on the internal homs by
\begin{align*}
\Map(\{c_i\},\{d_j\}) &= \lim_i \colim_j \Map(c_i, d_j) \to \lim_i \colim_j \Map(F(c_i), F(d_j)) \\
&\to \Map(\colim_i F(c_i), \colim_i F(d_j)).
\end{align*}
We saw in the preliminaries that $\wt F$ preserves all colimits and has a right adjoint (as functor of unenriched categories). In particular, it is part of an enriched adjunction if and only if it preserves tensors. To see that this is the case, let $X = \colim_j X_j$ be a simplicial set written as a filtered colimit of finite simplicial sets. Then $\{c_i\}_i \otimes X \cong \{c_i \otimes X_j\}_{(i,j)}$, hence $F(\{c_i\}_i \otimes X) \cong \colim_{(i,j)} F(c_i) \otimes X_j \cong \wt F(\{c_i\}) \otimes X$, using the hypothesis that $F$ preserves tensors with finite simplicial sets.
\end{proof}

In this section we give some assumptions under which an adjunction of the type above is a Quillen adjunction, and give a further criterion for this adjunction to be a Quillen equivalence. This gives a straightforward way of constructing ``profinite'' versions of certain classical Quillen adjunctions, as illustrated in \autoref{example:ProfiniteNerveQuillenAdjunction}. At the end of this section, we show that if $\bfC \subset \bfE$ inherits the structure of a (co)fibration test category in the sense of \autoref{example:InheritedCofibrationTestCategory}, then the ind- or pro-completion functor (relative to $\bfC$) is a Quillen functor.

\begin{definition}
A \emph{morphism of cofibration test categories} $\phi \colon (\bfC_1, \bfT_1) \to (\bfC_2, \bfT_2)$ is a simplicial functor $\phi \colon \bfC_1 \to \bfC_2$ that preserves finite colimits, finite tensors and (trivial) cofibrations, and in particular maps the full subcategory $\bfT_1$ into $\bfT_2$. Dually, a \emph{morphism of fibration test categories} $\phi \colon (\bfC_1, \bfT_1) \to (\bfC_2, \bfT_2)$ is a simplicial functor $\phi \colon \bfC_1 \to \bfC_2$ that preserves finite limits, finite cotensors and (trivial) fibrations and in particular maps the full subcategory $\bfT_1$ into $\bfT_2$.
\end{definition}

\begin{example}\label{example:NerveIsMorphismOfFibrationTestCats}
The nerve functor $N \colon \Fin\Grpd \to \LsSet_{KQ}$ is a morphism of fibration test categories. Similarly, taking the nerve of a category gives a morphism of fibration test categories $N \colon \Fin\Cat \to \LsSet_J$.
\end{example}

\begin{remark}\label{remark:RightAdjointPreservesFilteredColimits}
	If $\phi \colon (\bfC_1,\bfT_1) \to (\bfC_2,\bfT_2)$ is a morphism of cofibration test categories, then its canonical filtered colimit preserving extension $\phi_! \colon \Ind(\bfC_1) \to \Ind(\bfC_2)$ has a right adjoint $\phi^* \colon \Ind(\bfC_2) \to \Ind(\bfC_1)$ by \autoref{lem:TensorPreservingHasEnrichedRightAdjoint}. Since $\phi_!$ is an extension of $\phi \colon \bfC_1 \to \bfC_2$, it sends all objects in the image of $\bfC_1 \hookrightarrow \Ind(\bfC_1)$ to compact objects in $\Ind(\bfC_2)$, hence its right adjoint $\phi^*$ must preserve filtered colimits. Dually, if $\phi$ is a morphism of fibration test categories, then it canonically extends to a functor $\phi_* \colon \Pro(\bfC_1) \to \Pro(\bfC_2)$ that admits a right adjoint $\phi^* \colon \Pro(\bfC_2) \to \Pro(\bfC_1)$ which preserves cofiltered limits.
\end{remark}

\begin{proposition}\label{proposition:InducedQuillenPairs}
Let $\phi \colon (\bfC_1, \bfT_1) \to (\bfC_2, \bfT_2)$ be a morphism of cofibration test categories. Then the induced adjunction from \autoref{remark:RightAdjointPreservesFilteredColimits}
\[\phi_! : \Ind(\bfC_1) \rightleftarrows \Ind(\bfC_2) : \phi^*\]
is a simplicial Quillen adjunction. Dually, for a morphism of fibration test categories $\phi \colon (\bfC_1, \bfT_1) \to (\bfC_2, \bfT_2)$, the induced adjunction from \autoref{remark:RightAdjointPreservesFilteredColimits}
\[\phi^* : \Pro(\bfC_2) \rightleftarrows \Pro(\bfC_1) : \phi_*\]
is a simplicial Quillen adjunction.
\end{proposition}

\begin{proof}
Suppose $\phi \colon (\bfC_1, \bfT_1) \to (\bfC_2, \bfT_2)$ is a morphism of cofibration test categories. By \autoref{lem:TensorPreservingHasEnrichedRightAdjoint}, the adjunction $\phi_! \dashv \phi^*$ is an enriched adjunction of simplicial functors. Since $\phi_!$ extends $\phi$ and $\phi \colon \bfC_1 \to \bfC_2$ preserves all (trivial) cofibrations, we conclude that $\phi_! \colon \Ind(\bfC_1) \to \Ind(\bfC_2)$ preserves all generating (trivial) cofibrations. We conclude that $\phi_! \dashv \phi^*$ is a simplicial Quillen adjunction. The case of fibration test categories is dual.
\end{proof}

\begin{remark}
One could weaken the definition of a morphism of (co)fibration test categories $\phi \colon \bfC_1 \to \bfC_2$ by only asking it to be an (unenriched) functor of underlying categories and not asking it to preserve (co)tensors. In this case, one would still obtain a Quillen adjunction between the completed model structures, but it would merely be a Quillen adjunction between the underlying model categories, and not a simplicial one. Moreover, the proof of \autoref{proposition:InducedQuillenPairEquivalence} below would not go through in this case.
\end{remark}

\begin{example}\label{example:ProfiniteNerveQuillenAdjunction}
The nerve functors from \autoref{example:NerveIsMorphismOfFibrationTestCats} induce simplicial Quillen adjunctions
$\wh \Pi_1 : \s\wh\Set_Q \rightleftarrows \wh\Grpd : \wh N$
and
$\wh h : \s\wh\Set_J \rightleftarrows \wh\Cat : \wh N.$
These left adjoints are profinite versions of the fundamental groupoid and the homotopy category, respectively.
\end{example}

We call the restriction $\phi \colon \bfT_1 \to \bfT_2$ of a morphism of cofibration test categories \emph{homotopically essentially surjective} if for any $t' \in \bfT_2$, there exists a $t \in \bfT_1$ together with a weak equivalence $\phi(t) \wearrow t'$ in $\bfT_2$.

\begin{proposition}\label{proposition:InducedQuillenPairEquivalence}
Let $\phi \colon (\bfC_1, \bfT_1) \to (\bfC_2, \bfT_2)$ be a morphism of (co)fibration test categories.
\begin{enumerate}[(a)]
    \item \label{InducedQuillenPairEquivalence:DetectingWE} If the restriction $\bfT_1 \to \bfT_2$ of $\phi$ is homotopically essentially surjective, then $\phi^*$ detects weak equivalences.
    \item \label{InducedQuillenPairEquivalence:QuillenEquiv} In the case of a morphism of cofibration test categories, if moreover for any $t \in \bfT_1$ and $c \in \bfC_1$ the map 
    \[\Map(t,c) \to \Map(\phi(t), \phi(c)) \]
    is a weak equivalence, then the induced Quillen adjunction of \autoref{proposition:InducedQuillenPairs} is a Quillen equivalence.
    \item[(b')] \label{InducedQuillenPairEquivalence:QuillenEquivFibration} In the case of fibration test categories, if $\phi$ is homotopically essentially surjective and for any $t \in \bfT_1$ and $c \in \bfC_1$, the map
    \[\Map(c,t) \to \Map(\phi(c), \phi(t)) \]
    is a weak equivalence, then the induced Quillen adjunction of \autoref{proposition:InducedQuillenPairs} is a Quillen equivalence.
\end{enumerate}
\end{proposition}

\begin{proof}
We again only include a proof for cofibration test categories, as the case of a morphism of fibration test categories is dual. For item \eqref{InducedQuillenPairEquivalence:DetectingWE}, let $f \colon C \to D$ be a map in $\Ind(\bfC_2)$ and suppose that $\phi^*(f)$ is a weak equivalence in $\Ind(\bfC_1)$. If $t' \in \bfT_2$, then since $\phi \colon \bfT_1 \to \bfT_2$ is homotopically essentially surjective, there is a $t \in \bfT_1$ together with an equivalence $\phi(t) \wearrow t'$. Since $C$ and $D$ are fibrant in $\Ind(\bfC_2$), the map $\Map(t',C) \to \Map(t',D)$ is a weak equivalence if and only if $\Map(\phi(t),C) \to \Map(\phi(t),D)$ is so. Since $\phi_!$ extends $\phi$ and the adjunction $\phi_! \dashv \phi^*$ is enriched, we see that
\[\begin{tikzcd}
\Map(\phi(t),C) \ar[r] \ar[d,"\cong"] & \Map(\phi(t),D) \ar[d,"\cong"] \\
\Map(t, \phi^*(C)) \ar[r,"\sim"] & \Map(t,\phi^*(D))
\end{tikzcd}\]
commutes, hence $\Map(t',C) \to \Map(t',D)$ is a weak equivalence.

For item \eqref{InducedQuillenPairEquivalence:QuillenEquiv}, since the right adjoint $\phi^*$ detects weak equivalences by part \eqref{InducedQuillenPairEquivalence:DetectingWE}, it suffices to show that the unit $C \to \phi^* \phi_! C$ is a weak equivalence for every cofibrant $C$ in $\Ind(\bfC_1)$. Since $C$ is an cofiltered limit of objects in $\bfC_1$, by \autoref{remark:RightAdjointPreservesFilteredColimits} it is enough to show that $c \to \phi^* \phi_! c$ is a weak equivalence for every $c$ in $\bfC_1$. By definition of the weak equivalences in $\Ind(\bfC_1)$ and by the simplicial adjunction $\phi_! \dashv \phi^*$, this is equivalent to
\[\Map(t,c) \to \Map(\phi_!(t), \phi_!(c)) \cong \Map(\phi(t), \phi(c)) \]
being a weak equivalence, which holds by assumption.
\end{proof}

An interesting consequence of \autoref{proposition:InducedQuillenPairEquivalence} is that if, for a (co)fibration test category $(\bfC,\bfT)$, one ``enlarges'' $\bfC$ to a bigger category $\bfC'$ but keeps $\bfT$ the same, then one obtains Quillen equivalent model structures on $\Ind(\bfC)$ and $\Ind(\bfC')$ (or $\Pro(\bfC)$ and $\Pro(\bfC')$). The next example gives an illustration of this.

\begin{example}
Recall from \autoref{example:KanQuillenLeanFibrationTestCategory} that the category of lean simplicial sets $\LsSet$ inherits the structure of a fibration test category from $\s\Set_{KQ}$. We could give the category of degreewise finite simplicial sets $\s\Fin\Set$ a similar structure of a fibration test category, namely by defining the test objects to be the lean Kan complexes and the (trivial) fibrations to be those of $\LsSet_{KQ}$. That is, the test objects and the (trivial) fibrations of $\s\Fin\Set$ and of $\LsSet_{KQ}$ are identical. It is well known that the pro-categories $\Pro(\s\Fin\Set)$ and $\Pro(\LsSet) \simeq \s\wh\Set$ are not equivalent. However, the inclusion $\iota \colon \LsSet_{KQ} \hookrightarrow \s\Fin\Set$ is a morphism of fibration test categories that satisfies item \eqref{InducedQuillenPairEquivalence:QuillenEquivFibration} of \autoref{proposition:InducedQuillenPairEquivalence}, hence the induced adjunction
\[\iota^* : \Pro(\s\Fin\Set) \rightleftarrows \s\wh\Set_Q : \iota_*\]
is a Quillen equivalence.
\end{example}

The hypotheses for item \eqref{InducedQuillenPairEquivalence:QuillenEquiv} of \autoref{proposition:InducedQuillenPairEquivalence} can usually be weakened, namely if $\bfT$ is ``large enough'' in the following sense.

\begin{definition}\label{definition:ClosedUnderPullbackAlongCofib}
Let $(\bfC,\bfT)$ be a cofibration test category. We say that $\bfT$ is \emph{closed under pushouts along cofibrations} if, for any cofibration $r \cofarrow s$ in $\bfT$ and any map $r \to t$ in $\bfT$, the pushout $s \cup_r t$ is again contained in $\bfT$.

Dually, for a fibration test category $(\bfC,\bfT)$, we say that $\bfT$ is \emph{closed under pullbacks along fibrations} if, for any fibration $s \fibarrow r$ and any map $t \to r$ in $\bfT$, the pullback $s \times_r t$ is again contained in $\bfT$.
\end{definition}

This definition can be seen as ensuring that $\bfT$ has all finite homotopy (co)limits. If $\bfT$ is closed under pushouts along cofibrations, then it is enough to assume in item \eqref{InducedQuillenPairEquivalence:QuillenEquiv} that the restriction $\phi \colon \bfT_1 \to \bfT_2$ is \emph{homotopically fully faithful}, i.e.~that $\Map(s,t) \to \Map(\phi(s), \phi(t))$ is a weak equivalence for all $s,t \in \bfT_1$. The main ingredient is the following useful lemma.

\begin{lemma}\label{lemma:CofibrantObjectIsDirectedColimitTestObjects}
Let $(\bfC, \bfT)$ be a cofibration test category and suppose that $\bfT$ is closed under pushouts along cofibrations. Then any cofibrant object in $\Ind(\bfC)$ is a filtered colimit of objects in $\bfT$.
\end{lemma}

\begin{proof}
The ``fat small object argument'' of \cite{Makkai2014FatSmallObject} shows that if $C$ in $\Ind(\bfC)$ is cofibrant, then it is a retract of a colimit $\colim_{i \in I} c_i$ indexed by a directed poset $I$ that has a least element $\bot$, such that $c_\bot$ is the initial object $\varnothing$ and such that $c_\bot \to c_i$ is a (finite) composition of pushouts of generating cofibrations for any $i$. (This follows from Theorem 4.11 of \cite{Makkai2014FatSmallObject} together with the fact that all objects in $\bfT$ are compact.) In particular, since $\bfT$ is closed under pushouts along cofibrations, it follows that $c_i \in \bfT$ for every $i \in I$. Since ind-categories are idempotent complete, it follows that any retract of such a colimit is an object of $\Ind(\bfT)$ as well. In particular, any cofibrant object of $\Ind(\bfC)$ lies in $\Ind(\bfT)$.
\end{proof}

We leave it to the reader to dualize \autoref{lemma:CofibrantObjectIsDirectedColimitTestObjects} to the context of fibration test categories.

\begin{proposition}\label{proposition:HomotopicallyFullyFaithfulInducesQuillenEquiv}
Let $\phi \colon (\bfC_1, \bfT_1) \to (\bfC_2, \bfT_2)$ be a morphism of cofibration test categories (or fibration test categories) and suppose that $\bfT_1$ is closed under pushouts along cofibrations (or closed under pullbacks along fibrations, respectively). If the restriction $\phi \colon \bfT_1 \to \bfT_2$ is homotopically essentially surjective and homotopically fully faithful, then the induced Quillen adjunction of \autoref{proposition:InducedQuillenPairs} is a Quillen equivalence.
\end{proposition}

\begin{proof}
We prove the statement for ind-categories. As in the proof of \autoref{proposition:InducedQuillenPairEquivalence}, it suffices to show that the unit $C \to \phi^*\phi_! C$ is a weak equivalence for every cofibrant $C$ in $\Ind(\bfC_1)$. By \autoref{lemma:CofibrantObjectIsDirectedColimitTestObjects} any cofibrant object is a filtered colimit of objects of $\bfT_1$, so by \autoref{remark:RightAdjointPreservesFilteredColimits} it suffices to show that $t \to \phi^* \phi_! t$ is a weak equivalence for every $t \in \bfT_1$. This follows exactly as in the proof of \autoref{proposition:InducedQuillenPairEquivalence}.
\end{proof}

Recall from \autoref{ssec:ProAndIndCategories} that if $\bfE$ is a complete category and if $\bfC \subset \bfE$ is a small full subcategory closed under finite limits, then the functor $U \colon \Pro(\bfC) \to \bfE$ that sends a pro-object to its limit in $\bfE$ has a left adjoint $\wh{(\cdot)}_{\Pro}$, the pro-$\bfC$ completion functor. Dually, if $\bfE$ is cocomplete and $\bfC$ is closed under finite colimits, then the canonical functor $U \colon \Ind(\bfC) \to \bfE$ has a right adjoint $\wh{(\cdot)}_{\Ind}$. In the situation where $\bfE$ is a simplicial model category and $\bfC$ is a (co)fibration test category, these adjunctions are almost by definition Quillen pairs. Note that in the case of pro-categories, this is the Quillen pair mentioned in item \eqref{thm:ApproximateTheoremItem3} of \autoref{thm:ApproximateTheorem}.

\begin{proposition}\label{proposition:CompletionIsQuillenAdjunction}
Let $\bfE$ be a simplicial model category in which every object is fibrant and $\bfC \subset \bfE$ a full subcategory closed under finite colimits and finite tensors with the inherited structure of a cofibration test category (in the sense of \autoref{example:InheritedCofibrationTestCategory}). Then
\[U : \Ind(\bfC) \rightleftarrows \bfE : \wh{(\cdot)}_{\Ind} \]
is a simplicial Quillen adjunction. Dually, if every object in $\bfE$ is cofibrant and $\bfC \subset \bfE$ is a full subcategory closed under finite limits and finite cotensors, given the inherited structure of a fibration test category (as in \autoref{example:InheritedFibrationTestCategory}), then 
\[\wh{(\cdot)}_{\Pro} : \bfE \rightleftarrows \Pro(\bfC) : U  \]
is a simplicial Quillen adjunction.
\end{proposition}

\begin{proof}
The first adjunction arises by applying \autoref{lem:TensorPreservingHasEnrichedRightAdjoint} to the inclusion $\bfC \hookrightarrow \bfE$. We need to show that the left adjoint $U$ preserves the generating (trivial) cofibrations. Note that $U$ agrees with the inclusion $\bfC \hookrightarrow \bfE$ when restricted to $\bfC \subset \Ind(\bfC)$. Since the generating (trivial) cofibrations are defined as the (trivial) cofibrations in $\bfC \subset \bfE$ between cofibrant objects, they are preserved by $U$.

The case for pro-$\bfC$ completion follows dually.
\end{proof}

\begin{example}
The proposition above shows that the profinite completion functors for $\s\Set_{KQ}$ and $\Grpd$ are left Quillen. These Quillen adjunctions fit into a commutative diagram 
\begin{equation*} \label{diag:SeveralQuillenPairs}
\begin{tikzcd}[column sep = huge, row sep = 4.5 em]
\s\Set_{KQ} \ar[r,bend left = 20,"\Pi_1",""{name=A, below}] \ar[d,bend right,"\wh{(\cdot)}_{\Pro}"{left},""{name=C, left}] & \Grpd \ar[l,bend left = 20,"N",""{name=B,above}] \ar[from=A, to=B, symbol=\dashv] \ar[d,bend right = 20,"\wh{(\cdot)}_{\Pro}"{left},""{name=E, left}] \\
\s\wh\Set_Q \ar[r,bend left = 20,"\wh\Pi_1",""{name=A, below}] \ar[u,bend right = 20,"U"{right},""{name=D,right}] \ar[from=C, to=D, symbol=\dashv]  & \wh \Grpd \ar[l,bend left = 20, "\wh N",""{name=B,above}] \ar[from=A, to=B, symbol=\dashv] \ar[u,bend right ,"U"{right},""{name=F,right}] \ar[from=E, to=F, symbol=\dashv].
\end{tikzcd}
\end{equation*}
where $\wh N$ is the nerve adjunction from \autoref{example:ProfiniteNerveQuillenAdjunction}. There is a similar diagram of Quillen adjunctions for the (profinite) Joyal model structure and the model category of (profinite) categories.
\end{example}


\section{Bousfield localizations}\label{sec:BousfieldLocalizations}

Suppose we are given a cofibration test category $(\bfC, \bfT)$ and that we wish to shrink the full subcategory of test objects $\bfT$ to a smaller one $\bfT' \subset \bfT$. If $\bfT'$ is closed under finite pushout-products, then $(\bfC,\bfT')$ is a cofibration test category by \autoref{example:RestrictingToLessTestObjects}, hence we obtain two model structures $\Ind(\bfC,\bfT)$ and $\Ind(\bfC,\bfT')$ on the category $\Ind(\bfC)$. Since the (trivial) cofibrations of $(\bfC,\bfT')$ are those of $(\bfC,\bfT)$ between objects of $\bfT'$, the sets of generating (trivial) cofibrations of $\Ind(\bfC,\bfT')$ are contained in those of $\Ind(\bfC,\bfT)$. In particular, the identity functor is right Quillen when viewed as a functor $\Ind(\bfC,\bfT) \to \Ind(\bfC,\bfT')$. Since there are fewer weak equivalences in $\Ind(\bfC,\bfT)$ than in $\Ind(\bfC,\bfT')$, this right Quillen functor is close to being a right Bousfield localization. Recall that a right Bousfield localization of a model category is a model structure on the same category with the same class of fibrations, but with a larger class of weak equivalences. The model category $\Ind(\bfC,\bfT')$ is not necessarily a right Bousfield localization of $\Ind(\bfC,\bfT)$ since it has fewer generating trivial cofibrations, and hence it might have more fibrations than $\Ind(\bfC,\bfT)$. However, it is a general fact about model categories that in such a situation, there exists a model structure on $\Ind(\bfC)$ with the weak equivalences of $\Ind(\bfC,\bfT')$ and the fibrations of $\Ind(\bfC,\bfT)$.

\begin{lemma}\label{lemma:MixingModelStructures}
Let $\bfE_\alpha$ and $\bfE_\beta$ be cofibrantly generated model structures on the same category $\bfE$ and suppose that sets of generating cofibrations $I_\alpha$ and $I_\beta$ and sets of generating trivial cofibration $J_\alpha$ and $J_\beta$ respectively, are given. If $I_\alpha \subset I_\beta$ and $J_\alpha \subset J_\beta$, and if $\bfE_\alpha$ has more weak equivalences than $\bfE_\beta$, then there exists a cofibrantly generated model structure on $\bfE$ with the weak equivalences of $\bfE_\alpha$ and the fibrations of $\bfE_\beta$.
\end{lemma}

\begin{proof}
It easily follows by checking the hypotheses of Theorem 11.3.1 of \cite{Hirschhorn2003Model} that the sets $I_\alpha \cup J_\beta$ and $J_\beta$ determine a cofibrantly generated model structure on $\bfE$ in which the weak equivalences agree with those of $\bfE_\alpha$. This model structure has the desired properties. As an example, we check item (4b) of Theorem 11.3.1 of \cite{Hirschhorn2003Model}, and leave the other hypotheses to the reader. This comes down to showing that if $E \to F$ has the right lifting property with respect to $J_\beta$ and is a weak equivalence in $\bfE_\alpha$, then it must have the right lifting property with respect to $I_\alpha \cup J_\beta$. It suffices to show that $E \to F$ has the right lifting property with respect to $I_\alpha$. Since $E \to F$ has the right lifting property with respect to $J_\beta$, it has so with respect to $J_\alpha \subset J_\beta$, hence it is a fibration in $\bfE_\alpha$. Since it is also a weak equivalence in $\bfE_\alpha$, it follows that it has the right lifting property with respect to $I_\alpha$ and hence with respect to $I_\alpha \cup J_\beta$.
\end{proof}

If $\bfE$ is a simplicial model category with a given full subcategory $\bfT \subset \bfE$, then $R_{\bfT} \bfE$ denotes (if it exists) the right Bousfield localization of $\bfE$ in which a map $E \to E'$ is a weak equivalence if and only if $\Map(t,E) \to \Map(t,E')$ is a weak equivalence for every $t \in \bfT$. We call such a map a $\bfT$-colocal weak equivalence. Dually, $L_{\bfT} \bfE$ denotes (if it exists) the left Bousfield localization of $\bfE$ in which $E \to E'$ is a weak equivalence if and only if $\Map(E',t) \to \Map(E,t)$ is a weak equivalence for every $t \in \bfT$. Such a map is called a $\bfT$-local weak equivalence.

\begin{proposition}\label{proposition:BousfieldLocalizations}
Let $(\bfC,\bfT)$ be a cofibration test category and let $\bfT' \subset \bfT$ be a full subcategory. Then the right Bousfield localization $R_{\bfT'} \Ind(\bfC)$ exists and is cofibrantly generated.

Dually, if $(\bfC,\bfT)$ is a fibration test category and $\bfT' \subset \bfT$ a full subcategory, then the left Bousfield localization $L_{\bfT'} \Pro(\bfC)$ exists and is fibrantly generated.
\end{proposition}

\begin{proof}
We first prove the proposition in the special case that $\bfT'$ is closed under finite pushout-products, and then deduce the general case from this. In this special case, $(\bfC,\bfT')$ is a cofibration test category as in \autoref{example:RestrictingToLessTestObjects}, so we obtain a cofibrantly generated model category on $\Ind(\bfC,\bfT')$ in which the weak equivalences are the $\bfT'$-colocal ones. We also have the model structure on $\Ind(\bfC)$ corresponding to the cofibration test category $(\bfC,\bfT)$, which by construction has more generating (trivial) cofibrations than $\Ind(\bfC,\bfT')$. By applying \autoref{lemma:MixingModelStructures}, we obtain the desired right Bousfield localization $R_{\bfT'} \Ind(\bfC)$.

Now suppose that $\bfT'$ is not necessarily closed under finite pushout-products. Let $\bfT''$ be the smallest full subcategory of $\bfT$ that contains $\bfT'$ and that is closed under finite pushout-products and isomorphisms. This category can be obtained by repeatedly enlarging $\bfT'$ by adding all objects isomorphic to an object of the form $t' \otimes U \cup_{s' \otimes U} s' \otimes V$ to $\bfT'$, for $s' \cofarrow t'$ a cofibration in $\bfT'$ and $U \cofarrow V$ a cofibration of finite simplicial sets. This produces a sequence of full subcategories $\bfT' \subset \bfT'_1 \subset \bfT'_2 \subset \ldots \subset \bfT$ such that $\bfT'' = \cup_{n \in \bbN} \bfT'_n$. We claim that the $\bfT'$-colocal weak equivalences and the $\bfT''$-colocal weak equivalences in $\Ind(\bfC)$ agree. By the above inductive construction of $\bfT''$, it suffices to show that for any cofibration $s' \cofarrow t'$ of $(\bfC, \bfT)$ with $s', t' \in \bfT'$ and any cofibration $U \cofarrow V$ in $\FinsSet$, the map
\[\Map(t' \otimes U \cup_{s' \otimes U} s' \otimes V,C) \to \Map(t' \otimes U \cup_{s' \otimes U} s' \otimes V,D) \]
is a weak equivalence for any $\bfT'$-colocal weak equivalence $C \to D$. We leave this as an exercise to the reader, noting that these pushouts can be taken out of the mapping spaces to obtain homotopy pullbacks.
\end{proof}

\begin{example}
Let $\LsSet_{KQ}$ be the category of lean simplicial sets with the structure of a fibration test category as in \autoref{example:KanQuillenLeanFibrationTestCategory}. The model structure $\Pro(\LsSet_{KQ})$ then coincides with Quick's model structure $\s\wh\Set_{Q}$ under the equivalence of categories $\Pro(\LsSet) \simeq \s\wh\Set_Q$ by \autoref{corollary:QuicksModelStructureAgrees}. In particular, by \autoref{proposition:BousfieldLocalizations}, the left Bousfield localization $L_{\bfT} \s\wh\Set_Q$ exists for any collection of lean Kan complexes $\bfT$. If one takes $\bfT$ to consist of the spaces $K(\mathbb{F}_p,n)$ for all $n \in \bbN$, then one obtains a model structure on $\s\wh\Set$ in which the weak equivalences are the maps that induce equivalences in $\mathbb{F}_p$-cohomology and in which the cofibrations are the monomorphisms. This is exactly Morel's model structure on $\s\wh\Set$ for pro-$p$ spaces \cite{Morel1996EnsemblesProfinis}. In particular, this is an alternative to the construction in \autoref{example:MorelsModelStructures}.
\end{example}

\begin{example}\label{example:CompleteSegalAsBousfieldLoc}
Recall the Reedy model structure (with respect to Quick's model structure on $\s\wh\Set$) on $\biswhSet$ from \autoref{example:ReedyFibrationTestCategory}. By \autoref{proposition:ReedyModelStructuresCoincide}, this model structure can be obtained by applying \autoref{theorem:ModelStructureProC} to a certain fibration test category $\LsSet^{(2)}_R$. In particular, \autoref{proposition:BousfieldLocalizations} ensures that the left Bousfield localization $L_{\bfT} \biswhSet$ exists for any collection $\bfT$ of Reedy fibrant doubly lean simplicial sets. For example, one can take $\bfT$ to be the collection of all doubly lean bisimplicial sets that are complete Segal spaces in the sense of \cite{Rezk2001HomotopyTheoryOfHomotopyTheory}. This model structure will be called the \emph{model structure for complete Segal profinite spaces} and denoted $\biswhSet_{CSS}$. We will study this model structure in detail in \autoref{sec:CompleteSegalProfiniteVSProfiniteQuasiCats}. In particular, we will show in \autoref{proposition:QuillenEquivalencesCompleteSegalvsQuasiCat} that $\biswhSet_{CSS}$ is equivalent to the model structure for profinite quasi-categories $\s\wh\Set_J$ from \autoref{example:JoyalFibrationTestCategory}.
\end{example}

Note that \autoref{proposition:BousfieldLocalizations} was the last missing piece in the proof of \autoref{thm:ApproximateTheorem} (except for item \eqref{thm:ApproximateTheoremItem4} of that theorem, which follows from \autoref{thm:UnderlyingInftyCategoryPro}).

\begin{proof}[Proof of \autoref{thm:ApproximateTheorem}]
Let $\bfE$ be a simplicial model category in which every object is cofibrant and let $\bfC \subset \bfE$ be a small full subcategory of $\bfE$ closed under finite limits and cotensors by finite simplicial sets. Then $(\bfC, \bfT')$, where $\bfT' \subset \bfC$ is the full subcategory on the fibrant objects, inherits the structure of a fibration test category from $\bfE$ in the sense of \autoref{example:InheritedFibrationTestCategory}.

Now suppose $\bfT$ is any collection of fibrant objects in $\bfC$. By applying \autoref{theorem:ModelStructureProC} to $(\bfC,\bfT')$ and then applying \autoref{proposition:BousfieldLocalizations} (with $\bfT$ and $\bfT'$ interchanged), we obtain a model structure on $\Pro(\bfC)$ together with a (fibrantly generated) left Bousfield localization $\Pro(\bfC) \rightleftarrows L_{\bfT} \Pro(\bfC)$. The weak equivalences of $L_{\bfT} \Pro(\bfC)$ are by definition the $\bfT$-local equivalences. By \autoref{theorem:ModelStructureProC}, any object in $\Pro(\bfC)$ (and hence in $L_\bfT \Pro(\bfC)$) is cofibrant. By \autoref{proposition:CompletionIsQuillenAdjunction}, we obtain a simplicial Quillen adjunction $\bfE \rightleftarrows \Pro(\bfC)$ and hence a simplicial Quillen adjunction $\bfE \rightleftarrows L_\bfT \Pro(\bfC)$. We conclude that the model structure $L_\bfT \Pro(\bfC)$ satisfies items \eqref{thm:ApproximateTheoremItem1}-\eqref{thm:ApproximateTheoremItem3} of \autoref{thm:ApproximateTheorem}.
\end{proof}


\section{Example: complete Segal profinite spaces vs profinite quasi-categories}\label{sec:CompleteSegalProfiniteVSProfiniteQuasiCats}

Recall that in \autoref{example:JoyalFibrationTestCategory}, we defined the profinite Joyal model structure. In this section, we will define another candidate for the homotopy theory of profinite $\infty$-categories, namely a profinite version of Rezk's model category of complete Segal spaces. We then show that there are two Quillen equivalences between the model category of complete Segal profinite spaces and the profinite Joyal model structure. After establishing these Quillen equivalences, we characterize in both these model categories the weak equivalences between the fibrant objects as the essentially surjective and fully faithful maps, where fully faithfulness is defined in terms of the Quick model structure. It is worth mentioning that in \autoref{remark:UnderlyingInftyCategoryProfiniteInftyCats}, we moreover give a precise description of the underlying $\infty$-category of these model categories.

Let us start with a short review of the theory of complete Segal spaces, originally defined by Rezk in \cite{Rezk2001HomotopyTheoryOfHomotopyTheory}. Consider the category $\bisSet = \s\Set^{\Delta^{op}}$ of \emph{bisimplicial sets}, or \emph{simplicial spaces}, equipped with the Reedy model structure (with respect to the Kan-Quillen model structure on $\s\Set$). We denote this model category by $\bisSet_R$. Objects of $\bisSet$ have two simplicial parameters. We denote the ``inner'' one by $n,m, \ldots$ and refer to it as the \emph{space parameter}, and we denote the ``outer'' one (corresponding to the $\Delta^{op}$ in $\s\Set^{\Delta^{op}}$) by $s,t,r,\ldots$. For any pair of simplicial sets $X$ and $Y$, one can define the \emph{external product} $X \ultimes Y$ by $(X \ultimes Y)_{t,n} = X_t \times Y_n$. Note that the external product $\Delta^t \ultimes \Delta^n$ is the functor $\Delta^{op} \times \Delta^{op} \to \Set$ represented by $([t],[n])$. In particular, the internal hom of $\bisSet$ can be defined by $(Y^X)_{t,n} = \Hom((\Delta^t \ultimes \Delta^n) \times X, Y)$. This internal hom allows one to regard $\bisSet$ as a simplicial category in multiple ways; the two simplicial enrichments that we will use are given by
\[ \Map_1(X,Y) := (Y^X)_{\bullet,0} \quad \text{and} \quad \Map_2(X,Y) := (Y^X)_{0,\bullet}. \]
The category $\bisSet$ is tensored and cotensored with respect to both of these enrichments.

As described in \cite[\S 10 \& \S 12]{Rezk2001HomotopyTheoryOfHomotopyTheory}, one can localize the Reedy model structure on $\bisSet$ by the Segal maps
\[\Sp \Delta^t \ultimes \Delta^0 \cofarrow \Delta^t \ultimes \Delta^0, \]
where $\Sp \Delta^t = \Delta[0,1] \cup \ldots \cup \Delta[t-1,t]$ is the spine of the $t$-simplex. This gives the model category $\bisSet_{SS}$ for \emph{Segal spaces}. Localizing one step further, by the map
\[\{0\} \ultimes \Delta^0 \cofarrow J \ultimes \Delta^0, \]
gives the model category $\bisSet_{CSS}$ for \emph{complete Segal spaces}. Here $J$ is the nerve of the groupoid with two objects and exactly one isomorphism between any ordered pair of objects. It is part of a cosimplicial object $J^\bullet$ in $\s\Set$, $J^t$ being the nerve of the groupoid with $t+1$ objects and exactly one isomorphism between any ordered pair of objects.

All three of the model structures $\bisSet_R$, $\bisSet_{SS}$ and $\bisSet_{CSS}$ are $\s\Set_{KQ}$-enriched model structures with respect to the enrichment $\Map_2$ mentioned above.

The model category $\bisSet_{CSS}$ is Quillen equivalent to $\s\Set_J$. In fact, there are Quillen pairs in both directions, whose right Quillen functors are the evaluation at the inner coordinate $n=0$
\[ \ev_0 \colon \bisSet \to \s\Set; \quad (\ev_0 X)_t = X_{t,0} \]
and the singular complex functor with respect to $J^\bullet$
\[ \Sing \colon \s\Set \to \bisSet; \quad \Sing(X)_{t,n} = \Map(J^n,X)_t = \Hom(\Delta^t \times J^n, X). \]
These Quillen equivalences are described in detail in \cite{JoyalTierney2007QuasivsSegal}. One can prove, using the Quillen equivalence $\ev_0$ together with the fact that $\bisSet_{CSS}$ is a cartesian closed model category, that $\bisSet_{CSS}$ is a $\s\Set_J$-enriched model category with respect to the simplicial enrichment $\Map_1$ mentioned above. Both of the above right Quillen functors are simplicial functors that preserve cotensors with respect to this simplicial enrichment. This is explained in detail in the proof of Proposition E.2.2 of \cite{RiehlVerity2022Elements}.

Now let $\LsSet^{(2)}$ be the category of \emph{doubly lean} bisimplicial sets, i.e.~those bisimplicial sets $X$ for which $X_{t,n}$ is finite for each $t$ and $n$, and such that $X \cong \cosk_{t,n}(X)$ for some $t$ and $n$. Here $\cosk_{t,n} \colon \bisSet \to \bisSet$ is the functor that restricts $X \in \bisSet$ to a functor $\Delta^{op}_{\leq t} \times \Delta^{op}_{\leq n} \to \Set$ and then right Kan extends along $\Delta^{op}_{\leq t} \times \Delta^{op}_{\leq n} \hookrightarrow \Delta^{op} \times \Delta^{op}$. This agrees with the notion of doubly lean as defined at the end of \autoref{ssec:ProAndIndCategories}, and it follows from (the dual of) \autoref{thm:IndValuedPresheavesIsIndSkeletalPresheaves} that the inclusion $\LsSet^{(2)} \hookrightarrow \biswhSet$ extends to an equivalence $\Pro(\LsSet^{(2)}) \simeq \biswhSet$.

Each of the three model structures $\bisSet_R$, $\bisSet_{SS}$ and $\bisSet_{CSS}$ gives rise to the structure of a fibration test category on $\LsSet^{(2)}$ by the general scheme of \autoref{example:InheritedFibrationTestCategory}. We will mainly be interested in the Reedy and the complete Segal model structures, so denote the corresponding fibration test categories by $\LsSet^{(2)}_R$ and $\LsSet^{(2)}_{CSS}$ respectively.

\begin{definition}
The model structures on $\biswhSet$ obtained by applying \autoref{theorem:ModelStructureProC} to the fibration test categories $\LsSet^{(2)}_R$ and $\LsSet^{(2)}_{CSS}$ will be called the \emph{Reedy model structure for profinite spaces} and \emph{model structure for complete Segal profinite spaces}, and denoted $\biswhSet_R$ and $\biswhSet_{CSS}$, respectively. A fibrant object in $\biswhSet_{CSS}$ will be called a \emph{complete Segal profinite space}.
\end{definition}

Since we can view $\bisSet_{CSS}$ as a simplicial model category in two ways, the full subcategory $\LsSet^{(2)}_{CSS}$ inherits two different structures of a fibration test category, namely one with respect to the enrichment $\Map_1$ and one with respect to $\Map_2$. The (trivial) fibrations of both fibration test category structures agree, so they will induce the same model structures on $\Pro(\LsSet^{(2)}) \cong \biswhSet$. This shows that we can view $\biswhSet_{CSS}$ as a $\s\Set_J$-enriched model category through the enrichment $\Map_1$, and as a $\s\Set_{KQ}$-enriched model category through $\Map_2$.\footnote{In fact, one can show that $\biswhSet_{CSS}$ is a $\bisSet_{CSS}$-enriched model category, strengthening this statement.} In what follows, we will consider the simplicial enrichment $\Map_1$, since this one is compatible with the right Quillen functors $\ev_0$ and $\Sing$ discussed above.

By \autoref{proposition:CompletionIsQuillenAdjunction}, the profinite completion functor $\bisSet \to \biswhSet$ is a left Quillen functor, whose right adjoint is given by the functor $U \colon \biswhSet \to \bisSet$ that sends a bisimplicial profinite set to its underlying bisimplicial set. Levelwise, this is the functor that sends a profinite set to its underlying set.

Since $\LsSet^{(2)}_{CSS}$ has fewer test objects than $\LsSet^{(2)}_R$, we see that $\biswhSet_{CSS}$ has more weak equivalences than $\biswhSet_R$. By \autoref{appendix:proposition:CofibrationsAreMonosBisimplicial}, the cofibrations are the monomorphisms in both model structures, hence $\biswhSet_{CSS}$ is a left Bousfield localization of $\biswhSet_R$. In particular, the construction of the model structure $\biswhSet_{CSS}$ given in \autoref{example:CompleteSegalAsBousfieldLoc} agrees with the one given here.

The right Quillen functors $\ev_0$ and $\Sing$ mentioned above restrict to morphisms of fibration test categories between $\LsSet_J$ and $\LsSet^{(2)}_{CSS}$, where $\LsSet_J$ is the category of lean simplicial sets (with the fibration test category structure from \autoref{example:JoyalFibrationTestCategory}). This amounts to showing that $\ev_0$ maps doubly lean bisimplicial sets to lean simplicial sets, and that $\Sing$ maps lean simplicial sets to doubly lean bisimplicial sets. In the case of $\ev_0$, this follows directly from the definition, while the case of $\Sing$ requires some work.

\begin{lemma}
The functor $\Sing \colon \s\Set \to \bisSet$ takes lean simplicial sets to doubly lean bisimplicial sets.
\end{lemma}

\begin{proof}
Let $X$ be a lean simplicial set and suppose that $X$ is $n$-coskeletal. It suffices to show that $\Sing(X)_{\bullet,m}$ and $\Sing(X)_{t,\bullet}$ are both $n$-coskeletal and degreewise finite simplicial sets for any $t,m \in \bbN$. Since $J^m$ is a degreewise finite simplicial set for every $m$, we see that $\Sing(X)_{\bullet,m} = \Map(J^m,X)$ is an $n$-coskeletal degreewise finite simplicial set for every $m$. This automatically shows that $\Sing(X)_{t,\bullet}$ is a degreewise finite simplicial set as well. It therefore remains to show that, for every $n$-coskeletal simplicial set $X$ and every $t$, the simplicial set $\Sing(X)_{t,\bullet} \cong \Hom(J^\bullet \times \Delta^t, X) \cong \Hom(J^\bullet, X^{\Delta^t})$ is $n$-coskeletal. Since any cotensor $X^Y$ of an $n$-coskeletal simplicial set $X$ is again $n$-coskeletal, it suffices to prove the case $t=0$. To this end, let $\partial J^{k+1}$ denote the simplicial subset
\[\partial J^{k+1} = \bigcup_{x \in (\Delta^{k+1})_k} J^{k} \subset J^{k+1},\]
or equivalently, the left Kan extension of $J^\bullet \colon \Delta \to \s\Set$ along the Yoneda embedding $\Delta \to \s\Set$, evaluated at $\partial \Delta^{k+1} \in \s\Set$. The inclusion $\partial J^{k+1} \hookrightarrow J^{k+1}$ restricts to an isomorphism $\sk_n \partial J^{k+1} \to \sk_n J^{k+1}$ for any $k \geq n$. Combining this with the canonical isomorphism $\Hom(\partial \Delta^{k+1}, \Hom(J^\bullet,X)) \cong \Hom(\partial J^{k+1}, X)$, it follows that $\Hom(J^\bullet,X)$ is $n$-coskeletal.
\end{proof}

Denote the profinite Joyal model structure by $\s\wh\Set_J$. We can apply \autoref{proposition:HomotopicallyFullyFaithfulInducesQuillenEquiv} to $\ev_0 \colon \LsSet^{(2)}_{CSS} \to \LsSet_J$ and $\Sing \colon \LsSet_J \to \LsSet^{(2)}_{CSS}$ to show that the induced functors between $\s\wh\Set_J$ and $\biswhSet_{CSS}$ are right Quillen equivalences. We will denote these functors by $\ev_0$ and $\Sing$ as well.

\begin{proposition}\label{proposition:QuillenEquivalencesCompleteSegalvsQuasiCat}
The functors $\ev_0 \colon \biswhSet_{CSS} \to \s\wh\Set_J$ and $\Sing \colon \s\wh\Set_J \to \biswhSet_{CSS}$ are right Quillen equivalences.
\end{proposition}

\begin{proof}
Since there is a natural isomorphism $\ev_0 \Sing(X) \cong X$, it suffices to show that $\ev_0 \colon \biswhSet_{CSS} \to \s\wh\Set_J$ is a right Quillen equivalence. The same then follows for $\Sing$ by the 2 out of 3 property. Since $\ev_0 \colon \bisSet_{CSS} \to \s\Set_J$ is a (simplicial) right Quillen equivalence, its restriction $\ev_0 \colon \LsSet^{(2)}_{CSS} \to \LsSet_J$ is a morphism of fibration test categories that is homotopically fully faithful when restricted to test objects. Furthermore, it is homotopically essentially surjective since $X \cong \ev_0(\Sing X)$ for any lean quasi-category $X$. By \autoref{proposition:HomotopicallyFullyFaithfulInducesQuillenEquiv}, we conclude that induced functor $\ev_0 \colon \biswhSet_{CSS} \to \s\wh\Set_J$ (and hence $\Sing \colon \s\wh\Set_J \to \biswhSet_{CSS}$) is a right Quillen equivalence.
\end{proof}

One can prove ``profinite versions'' of many of the properties that complete Segal spaces enjoy. The general strategy for proving such a profinite version of a given property is to reduce it to its classical counterpart. We will illustrate this by showing that the weak equivalences between complete Segal profinite spaces coincide with (a profinite version of) the Dwyer-Kan equivalences. This is done by exploiting two facts: that $\biswhSet_{CSS}$ is a left Bousfield localization of the Reedy model structure $\biswhSet_R$ (with respect to $\s\wh\Set_Q$), and that the weak equivalences between fibrant objects in $\s\wh\Set_Q$ can be detected underlying in $\s\Set_{KQ}$. Denote the functor that sends a simplicial profinite set to its underlying simplicial set by $U \colon \s\wh\Set \to \s\Set$. Note that this functor is right Quillen as functor from Quick's model structure to the Kan-Quillen model structure, and that its left adjoint is the profinite completion functor.

\begin{proposition}\label{proposition:QuickEquivalenceDetectedUnderlying}
A map $X \to Y$ between fibrant objects in $\s\wh\Set_Q$ is a weak equivalence if and only if $UX \to UY$ is a weak equivalence in $\s\Set_{KQ}$.
\end{proposition}

\begin{proof}
This follows from Theorem E.3.1.6 of \cite{lurie2016spectral}, which states that the functor between the underlying $\infty$-categories of $\s\wh\Set_Q$ and $\s\Set_{KQ}$ induced by $U$ (which is called ``$\mathrm{Mat}$'' by Lurie) is conservative. Another way to deduce this proposition is to show that the weak equivalences between fibrant objects in $\s\wh\Set_Q$ are the $\pi_*$-isomorphisms (as in the proof of Proposition 3.9 of \cite{BoavidaHorelRobertson2019Operads}) and that the the underlying group/set $U \pi_n(X,x)$ of the profinite group/set $\pi_n(X,x)$ agrees with $\pi_n(UX,x)$ for any fibrant $X \in \s\wh\Set_Q$ and any $x \in X_0$.
\end{proof}

Since $\biswhSet_{CSS}$ is a left Bousfield localization of the model category $\biswhSet_R$, which coincides with the Reedy model structure on $\biswhSet$ with respect to $\s\wh\Set_Q$ by \autoref{proposition:ReedyModelStructuresCoincide}, we see that a map between complete Segal profinite spaces is a weak equivalence if and only if it is levelwise a weak equivalence in $\s\wh\Set_Q$. In particular, we obtain the following result:

\begin{proposition}\label{proposition:CompleteSegalEquivalenceDetectedUnderlying}
A map $X \to Y$ between complete Segal profinite spaces is a weak equivalence if and only if for every $t$, the map $X_{t,\bullet} \to Y_{t,\bullet}$ is a weak equivalence in $\s\wh\Set_Q$. In particular, $X \to Y$ is a weak equivalence between complete Segal profinite spaces if and only if $UX \to UY$ is a weak equivalence between complete Segal spaces.
\end{proposition}

For a complete Segal profinite space $X$ and two objects $x,y \in X_{0,0}$, i.e.~two maps $\Delta^0 \ultimes \Delta^0 \to X$, we can mimic the classical definition of the mapping space $\map_X(x,y)$ by defining $\map_X(x,y)$ as the pullback
\[\begin{tikzcd} \map_X(x,y) \arrow[dr,phantom, very near start, "\lrcorner"] \ar[r] \ar[d] & X_{1,\bullet} \ar[d,"(d_1{,}d_0)"] \\
\Delta^0 \ar[r,"(x{,}y)"] & X_{0,\bullet} \times X_{0,\bullet}.
\end{tikzcd}\]
Since $X$ is Reedy fibrant, the map $(d_1,d_0) \colon X_{1,\bullet} \to X_{0,\bullet} \times X_{0,\bullet}$ is a fibration in $\s\wh\Set_Q$, and hence $\map_X(x,y)$ is a fibrant object in $\s\wh\Set_Q$. Since $U \colon \biswhSet_{CSS} \to \bisSet_{CSS}$ preserves limits, we see that $U(\map_X(x,y)) \cong \map_{UX}(x,y)$ for any complete Segal profinite space $X$. If $f \colon X \to Y$ is a map between complete Segal profinite spaces, then for any $x,y \in X_{0,0}$, we obtain a map $\map_X(x,y) \to \map_Y(fx,fy)$ from the universal property of the pullback. We call a map between complete Segal profinite spaces $f \colon X \to Y$ \emph{fully faithful} if, for any $x,y \in X_{0,0}$, the map $\map_X(x,y) \to \map_Y(fx,fy)$ is a weak equivalence in $\s\wh\Set_Q$. It follows from \autoref{proposition:QuickEquivalenceDetectedUnderlying} that $X \to Y$ is fully faithful if and only if $UX \to UY$ is a fully faithful map of complete Segal spaces.

One can also mimic the classical definitions of a homotopy and of homotopy equivalences in a complete Segal space, and use this to define what it means for a map of complete Segal profinite spaces to be essentially surjective. An equivalent, but easier, way is to say that $X \to Y$ is \emph{essentially surjective} if and only if the induced map $\pi_0 X_{0,\bullet} \to \pi_0Y_{0,\bullet}$ is an epimorphism of profinite sets. Since $U \pi_0 Z \cong \pi_0 U Z$ for any fibrant object $Z$ in $\s\wh\Set_Q$, and since epimorphisms of profinite sets are detected underlying, we see that a map of complete Segal profinite spaces $X \to Y$ is essentially surjective if and only if $UX \to UY$ is so.

\begin{definition}\label{definition:DKEquivalenceSegal}
A map between complete Segal profinite spaces is called a \emph{Dwyer-Kan equivalence} or \emph{DK-equivalence} if it is essentially surjective and fully faithful.
\end{definition}

\begin{theorem}\label{theorem:DKEquivalenceSegal}
A map between complete Segal profinite spaces is a Dwyer-Kan equivalence if and only if it is a weak equivalence in $\biswhSet_{CSS}$.
\end{theorem}

\begin{proof}
As explained above \autoref{definition:DKEquivalenceSegal}, $f \colon X \to Y$ is essentially surjective and fully faithful if and only if $UX \to UY$ is so. By Proposition 7.6 of \cite{Rezk2001HomotopyTheoryOfHomotopyTheory}, this is the case if and only if $UX \to UY$ is a weak equivalence in $\bisSet_{CSS}$. By \autoref{proposition:CompleteSegalEquivalenceDetectedUnderlying}, this is equivalent to $X \to Y$ being a weak equivalence in $\biswhSet_{CSS}$.
\end{proof}

One can lift \autoref{proposition:CompleteSegalEquivalenceDetectedUnderlying} and \autoref{theorem:DKEquivalenceSegal} to analogous results about weak equivalences between profinite quasi-categories using the Quillen equivalences $\ev_0$ and $\Sing$ between $\s\wh\Set_J$ and $\biswhSet_{CSS}$.

\begin{proposition}\label{proposition:ProfiniteQuasiCatEquivalenceDetectedUnderlying}
A map $X \to Y$ between profinite quasi-categories is a weak equivalence in $\s\wh\Set_J$ if and only if $UX \to UY$ is a weak equivalence in $\s\Set_J$.
\end{proposition}

\begin{proof}
Let $f \colon X \to Y$ be a map between profinite quasi-categories. If $f$ is a weak equivalence in $\s\wh\Set_J$, then $Uf \colon UX \to UY$ is a weak equivalence of quasi-categories since $U \colon \s\wh\Set_J \to \s\Set_J$ is right Quillen. Conversely, suppose $Uf$ is a weak equivalence of quasi-categories. Then $\Sing Uf \colon \Sing(UX) \to \Sing(UY)$ is a weak equivalence between complete Segal spaces. Note that $\Sing \circ U \simeq U \circ \Sing$ since both functors preserve cofiltered limits and since they agree on lean simplicial sets. By \autoref{proposition:CompleteSegalEquivalenceDetectedUnderlying}, we see that $\Sing X \to \Sing Y$ is a weak equivalence between complete Segal profinite spaces. Since $\ev_0$ is right Quillen, we see that the original map $\ev_0 \Sing X \cong X \to Y \cong \ev_0 \Sing Y$ is a weak equivalence in $\s\wh\Set_J$.
\end{proof}

For a profinite quasi-category $X$ and two $0$-simplices $x,y \in X_0$ (i.e.~maps $\Delta^{0} \to X$), we define $\map_X(x,y)$ as the pullback
\[\begin{tikzcd}
\map_X(x,y) \arrow[dr,phantom, very near start, "\lrcorner"] \ar[r] \ar[d] & X^{\Delta^1} \ar[d,"(\ev_0{,}\ev_1)"] \\
\Delta^0 \ar[r,"(x{,}y)"] & X \times X
\end{tikzcd}\]
Since the right-hand vertical map is obtained by cotensoring with the cofibration $\partial \Delta^1 \hookrightarrow \Delta^1$, it must be a fibration in $\s\wh\Set_J$. In particular, $\map_X(x,y)$ is fibrant in $\s\wh\Set_J$. One can show that, analogous the classical case, $\map_X(x,y)$ is actually fibrant in $\s\wh\Set_Q$. However, the proof of this is technical and not necessary for what follows, so it is not included.

A map $f \colon X \to Y$ of profinite quasi-categories induces a morphism $\map_X(x,y) \to \map_Y(fx,fy)$ for any $x,y \in X_0$ by the universal property of the pullback. We say that $f$ is \emph{fully faithful} if $\map_X(x,y) \to \map_Y(fx,fy)$ is a weak equivalence in $\s\wh\Set_J$ for any $x,y \in X_0$.\footnote{Since the simplicial profinite sets $\map_X(x,y)$ and $\map_Y(fx,fy)$ are actually fibrant in $\s\wh\Set_Q$, this is equivalent to asking that $\map_X(x,y) \to \map_Y(fx,fy)$ is a weak equivalence in $\s\wh\Set_Q$.} For a $1$-simplex $\alpha \in X_1$ with $d_1 \alpha = x$ and $d_0 \alpha = y$, i.e a $0$-simplex in $\map_X(x,y)$, we say that $\alpha$ is a \emph{homotopy equivalence} if $\Delta^1 \xrightarrow{\alpha} X$ extends to a map $J^1 \to X$. Here $J^1$ is viewed as a simplicial profinite set through the inclusion $\s\Fin\Set \hookrightarrow \s\wh\Set$. We say that a map of profinite quasi-categories $f \colon X \to Y$ is \emph{essentially surjective} if for any $y \in Y_0$, there exists an $x \in X_0$ and an $\alpha \in \map_Y(fx,y)$ such that $\alpha$ is a homotopy equivalence.

Since $U \colon \s\wh\Set \to \s\Set$ preserves pullbacks, we see that $U \map_X(x,y) \cong \map_{UX}(x,y)$. By \autoref{proposition:ProfiniteQuasiCatEquivalenceDetectedUnderlying}, a map $X \to Y$ of profinite quasi-categories is fully faithful if and only if $UX \to UY$ is so. Since $\Hom(J^1,X) \cong \Hom(J^1,UX)$ for any $X \in \s\wh\Set$, we also see that $X \to Y$ is essentially surjective if and only if $UX \to UY$ is so.

\begin{definition}\label{definition:DKEquivalenceQuasi}
A map between profinite quasi-categories is called a \emph{Dwyer-Kan equivalence} or \emph{DK-equivalence} if it is essentially surjective and fully faithful.
\end{definition}

\begin{theorem}
A map between profinite quasi-categories is a Dwyer-Kan equivalence if and only if it is a weak equivalence in $\s\wh\Set_J$.
\end{theorem}

\begin{proof}
A map $X \to Y$ of profinite quasi-categories is a DK-equivalence if and only if $UX \to UY$ is. Since the weak equivalences between fibrant objects in $\s\Set_J$ are exactly the DK-equivalences, we conclude from \autoref{proposition:ProfiniteQuasiCatEquivalenceDetectedUnderlying} that a map of profinite quasi-categories $X \to Y$ is a DK-equivalence if and only if it is a weak equivalence.
\end{proof}

\appendix

\section{Comparison to the \texorpdfstring{$\infty$}{infinity}-categorical approach}\label{appendix:InftyApproach}

The goal of this appendix is to compare the model structures on $\Ind(\bfC)$ and $\Pro(\bfC)$ constructed in this paper to the $\infty$-categorical approach to ind- and pro-categories. Since the cases of ind- and pro-categories are dual, we only treat the case of ind-categories and dualize the main result at the end of this appendix.

Given a cofibration test category $\bfC$, the underlying $\infty$-category of the completed model structure on $\Ind(\bfC)$ will be denoted by $\Ind(\bfC)_\infty$. Recall that this $\infty$-category is defined as the homotopy-coherent nerve of the full simplicial subcategory spanned by the fibrant-cofibrant objects. We will show that if $(\bfC,\bfT)$ is a cofibration test category with a suitable assumption on $\bfT$, then the $\infty$-category $\Ind(\bfC)_\infty$ is equivalent to $\Ind(N(\bfT))$. Here $N(\bfT)$ is the homotopy-coherent nerve of the simplicial category $\bfT$, and $\Ind$ denotes the $\infty$-categorical version of the ind-completion as defined in Definition 5.3.5.1 of \cite{Lurie2009HTT}.

\begin{warrning}
There is a subtlety here that we should point out: if $(\bfC,\bfT)$ is a cofibration test category with respect to the Joyal model structure on $\s\Set$, meaning that items \eqref{CofTestCat:TestObjectsClosedUnderTensors} and \eqref{CofTestCat:TestingTrivialCofibration} of \autoref{definition:CofTestCat} hold with respect to the trivial cofibrations and weak equivalences of $\s\Set_J$, then the ``mapping spaces'' of $\bfT$ are quasi-categories but not necessarily Kan complexes. Recall that any quasi-category $X$ contains a maximal Kan complex, which we will denote by $k(X)$.
Since this functor $k$ preserves cartesian products, any category enriched in quasi-categories can be replaced by a category enriched in Kan complexes by applying the functor $k$ to the simplicial hom.
If $(\bfC,\bfT)$ is a cofibration test category with respect to $\s\Set_J$, then we will abusively write $N(\bfT)$ for the simplicial set obtained by first applying the functor $k$ to all the mapping spaces in $\bfT$, and then applying the homotopy-coherent nerve. Similarly, by the underlying infinity category $\Ind(\bfC)_\infty$ of $\Ind(\bfC)$, we mean the quasi-category obtained by taking the full subcategory on fibrant-cofibrant objects, applying $k$ to all mapping spaces, and then taking the homotopy-coherent nerve.
\end{warrning}

Since $\Ind(\bfC)_\infty$ is the underlying $\infty$-category of a combinatorial model category, we see that it is complete and cocomplete. Furthermore, since $\bfT$ is a full subcategory of the fibrant-cofibrant objects in $\Ind(\bfC)$, we see that the inclusion $\bfT \hookrightarrow \Ind(\bfC)$ induces a fully faithful inclusion $N(\bfT) \hookrightarrow \Ind(\bfC)_\infty$. By Proposition 5.3.5.10 of \cite{Lurie2009HTT}, this inclusion extends canonically to a filtered colimit preserving functor $F \colon \Ind(N(\bfT)) \to \Ind(\bfC)_\infty$. In order for this functor to be an equivalence, any object in $\Ind(\bfC)$ needs to be equivalent to a filtered homotopy colimit of objects in $\bfT$. This means that $\bfT$ should be ``large enough'' for this to hold. It turns out that this is the case if $\bfT$ is closed under pushouts along cofibrations (in the sense of \autoref{definition:ClosedUnderPullbackAlongCofib}).

\begin{theorem}\label{thm:UnderlyingInftyCategory}
Let $(\bfC,\bfT)$ be a cofibration test category and suppose that $\bfT$ is closed under pushouts along cofibrations. Then the canonical functor
\[F \colon \Ind(N(\bfT)) \to \Ind(\bfC)_\infty\]
is an equivalence of quasi-categories.
\end{theorem}

\begin{remark}
Note that in many of the examples discussed in this paper, the category $\bfT$ of test objects is closed under pushouts along cofibrations. For example, this is the case if $(\bfC,\bfT)$ has inherited the structure of a cofibration test category from some model category $\bfE$ in the sense of \autoref{example:InheritedCofibrationTestCategory}.
\end{remark}

\begin{remark}
If $(\bfC,\bfT)$ is a cofibration test category, then one can always ``enlarge'' the full subcategory $\bfT$ together with the sets of (trivial) cofibrations to obtain a cofibration test category $(\bfC,\bfT')$ such that $\bfT'$ is closed under pushouts along cofibrations, and for which the completed model structures $\Ind(\bfC,\bfT)$ and $\Ind(\bfC,\bfT')$ coincide. To see this, note that we can define $\bfT'$ to consists of all objects in $\bfC$ that are cofibrant in $\Ind(\bfC,\bfT)$, and that we can define the (trivial) cofibrations of $(\bfC,\bfT')$ to be the trivial cofibrations of $\Ind(\bfC,\bfT)$ between objects of $\bfT'$; that is, we endow $\bfC$ with the structure of a cofibration test category inherited from $\Ind(\bfC,\bfT)$ (see \autoref{example:InheritedCofibrationTestCategory}). It is then clear that the model structures $\Ind(\bfC,\bfT)$ and $\Ind(\bfC,\bfT')$ coincide, and that $\bfT'$ is closed under pushouts along cofibrations. In particular, we see by \autoref{thm:UnderlyingInftyCategory} that the underlying $\infty$-category of $\Ind(\bfC,\bfT)$ can be described as the ind-category of the small $\infty$-category $N(\bfT')$, which contains $N(\bfT)$ as a full subcategory.
\end{remark}

Before proving this theorem, we will prove the following rectification result.

\begin{lemma}\label{lem:Rectification}
Let $(\bfC,\bfT)$ be a cofibration test category such that $\bfT$ is closed under pushouts along cofibrations, and let $I$ be a poset with the property that $I_{< i}$ is finite for every $i$. For any diagram $X \colon N(I) \to N(\bfT)$, there exists a strict diagram $Y \colon I \to \bfT$ such that $N(Y) \colon N(I) \to N(\bfT)$ is naturally equivalent to $X$. This diagram $Y$ can be constructed such a way that or any $i \in I$, the map
\[\colim_{j < i} Y_j \to Y_i \]
is a composition of two pushouts of cofibrations in $\bfT$.
\end{lemma}

The following lemma is needed for the proof.

\begin{lemma}\label{lem:ContinuousColimIsGood}
Let $(\bfC,\bfT)$ be a cofibration test category and let $\{Y_i\}_{i \in I}$ be a diagram in $\bfT$ indexed by a finite poset such that for any $i \in I$, the map
\[\colim_{j < i} Y_j \to Y_i \]
is a finite composition of pushouts of cofibrations of $(\bfC,\bfT)$. Then, for any $k \in I$, the map
\[Y_k \to \colim_{i \in I} Y_i \]
is a finite composition of pushouts of cofibrations. In particular, if $\bfT$ is closed under pushouts along cofibrations, then $\colim_i Y_i$ is an object of $\bfT$.
\end{lemma}

\begin{proof}
This follows from the dual of \cite[Proposition 2.17]{BarneaSchlank2016ProSimplicialSheaves}. For the convenience of the reader, we spell out their argument in our setting. Throughout this proof, we call a map in $\bfC$ \emph{good} if it is a finite composition of pushouts of cofibrations. Note that any pushout of a good map is again a good map. A subposet $S \subset I$ is called a \emph{sieve} if for any $i \in S$ and any $j \leq i$ in $I$, one has $j \in S$. Write $Y_S = \colim_{j \in S} Y_j$ for any sieve $S$ and $Y_{<i}$ for $Y_{I_{<i}} = \colim_{j<i} Y_j$ for any $i \in I$.

We will prove inductively that for two sieves $S \subset T$, the map $Y_S \to Y_T$ is good. This certainly holds if $|T| = 0$, so suppose this holds for $|T| < n$ and let sieves $S \subset T$ with $|T| = n$ be given. If $S = T$ then there is nothing to prove, so suppose that $S \subsetneq T$ and choose some maximal $i \in T \setminus S$. We then obtain a diagram
\[\begin{tikzcd}
 & Y_{< i} \ar[r,"\mathrm{good}"] \ar[d] & Y_{i} \ar[d] \\
Y_S \ar[r,"\mathrm{good}"] & Y_{T \setminus \{i\}} \ar[r] & Y_{T}. \ar[ul, phantom, very near start, "\ulcorner"]
\end{tikzcd}\]
where the square is a pushout. The map $Y_{<i} \to Y_i$ is good by assumption while $Y_S \to Y_{T \setminus \{i\}}$ is good by the induction hypothesis, so we conclude that $Y_S \to Y_T$ is good. This completes the induction and the lemma now follows by considering the sieves $S = I_{\leq k}$ and $T = I$.
\end{proof}

\begin{proof}[Proof of \autoref{lem:Rectification}]
To distinguish colimits in quasi-categories from homotopy colimits and ordinary colimits in simplicial categories, we will call them $\infty$-colimits. By a homotopy colimit of a diagram $Z \colon J \to \bfT$, we mean a cocone $Z_j \to W$ that induces an equivalence
\[\Map(W,t) \wearrow \holim_{j \in J} \Map(Z_j,t) \]
for every $t \in \bfT$. The following proof is for the case that $(\bfC,\bfT)$ is a cofibration test category with respect to the Kan-Quillen model structure on $\s\Set$. The same proof works if $(\bfC,\bfT)$ is a cofibration test category with respect to $\s\Set_J$; however, one has to replace $\Map(-,-)$ with the maximal Kan complex $k(\Map(-,-))$ contained in it, and one has to replace $\Delta^1$ by the simplicial set $H$ (as defined in \autoref{lem:SkeletalJ}) in the construction of the mapping cylinder below.

We will construct the diagram $Y \colon I \to \bfT$ and the equivalence $N(Y) \simeq X$ inductively. Let $i \in I$ be given and suppose that $Y|_{I_{<i}} \colon I_{< i} \to \bfT$ and $N(Y|_{I_{<i}}) \simeq X|_{N(I_{<i})}$ have been constructed and have the desired properties. We need to construct $Y|_{I_{\leq i}} \colon I_{\leq i} \to \bfT$ and an equivalence $N(Y_{I_{\leq i}}) \simeq X|_{I_{\leq i}}$ extending these. Denote $\colim_{j < i} Y_j$ by $Y_{<i}$. If $I_{< i}$ is empty, then $Y_{<i}$ is the initial object of $\bfC$ and hence an object of $\bfT$ by definition. If $I_{<i}$ is not empty, then it follows from the assumptions on $Y|_{I_{<i}}$ and \autoref{lem:ContinuousColimIsGood} that $Y_{<i}$ is an object of $\bfT$. Also note that the assumptions on $Y|_{I_{<i}}$ and the fact that $\Ind(\bfC)$ is a simplicial model structure ensure that, for any $t \in \bfT$, the diagram $j \mapsto \Map(Y_j, t)$ is fibrant in the injective model structure on $\s\Set^{(I_{< i})^{op}}$. In particular, we see that
\[\Map(Y_{<i},t) \cong \lim_{j < i} \Map(Y_j,t) \simeq \holim_{j < i} \Map(Y_j, t),\]
so $Y_{<i}$ is a homotopy colimit of the diagram $Y|_{I_{<i}}$. By Theorem 4.2.4.1 of \cite{Lurie2009HTT}, it follows that it is also the $\infty$-colimit of the diagram $N(Y|_{I_{<i}}) \colon N(I_{<i}) \to N(\bfT)$. In particular, if we define the diagram $Y' \colon I_{\leq i} \to \bfT$ by $Y'_j = Y_j$ for all $j < i$ and $Y'_i = Y_{<i}$, then the natural equivalence $N(Y|_{I_{<i}}) \simeq X|_{N(I_{<i})}$ extends to a natural map $N(Y') \to X|_{N(I_{\leq i})}$. The map $Y_{<i} = Y'_i \to X_i$ factors through the mapping cylinder
\[Y_{<i} \cong Y_{<i} \otimes \{0\} \to Y_{<i} \otimes \Delta^1 \cup_{Y_{<i} \otimes \{1\}} X_i \otimes \{1\} \wearrow X_i \]
in $\bfT$, where the second map is a weak equivalence. The first map can be written as a composition of the following two pushouts of cofibrations:
\[\begin{tikzcd}
\varnothing \ar[d,tail] \ar[r] & Y_{<i} \otimes \{0\} \ar[d] & Y_{< i} \otimes \partial \Delta^1 \ar[d,tail] \ar[r] & Y_{<i} \otimes \{0\} \sqcup X_i \otimes \{1\} \ar[d] \\
X_i \otimes \{1\} \ar[r] & Y_{<i} \otimes \{0\} \sqcup X_i \otimes \{1\}, \ar[ul, phantom, very near start, "\ulcorner"] & Y_{< i} \otimes \Delta^1 \ar[r] &  Y_{<i} \otimes \Delta^1 \underset{Y_{<i} \otimes \{1\}}{\cup} X_i \otimes \{1\}. \ar[ul, phantom, very near start, "\ulcorner"]
\end{tikzcd} \]
Define $Y_i = Y_{<i} \otimes \Delta^1 \cup_{Y_{<i} \otimes \{1\}} X_i$. This defines a diagram $Y|_{I_{\leq i}} \colon I_{\leq i} \to \bfT$. The above factorization of $Y_{<i} \to X_i$ shows that we obtain a natural equivalence $N(Y|_{I_{\leq i}}) \simeq X|_{N(I_{\leq i})}$ extending the equivalence $N(Y|_{I_{< i}}) \simeq X|_{N(I_{< i})}$.
\end{proof}

We are now ready to prove \autoref{thm:UnderlyingInftyCategory}.

\begin{proof}[Proof of \autoref{thm:UnderlyingInftyCategory}]
The terms ``colimit'', ``homotopy colimit'' and ``$\infty$-colimit'' are used in the same way as in the proof of \autoref{lem:Rectification}. We will denote mapping spaces in a simplicial category by ``$\Map$'', while mapping spaces in a quasi-category are denoted by ``$\map$''; that is, with a lowercase m.

We will prove that the functor $F \colon \Ind(N(\bfT)) \to \Ind(\bfC)_\infty$ is fully faithful and essentially surjective. To see that $F$ is fully faithful, we need to show that
\[\map_{\Ind(N(\bfT))}(X,Y) \to \map_{\Ind(\bfC)_\infty}(F(X),F(Y)) \]
is a weak equivalence for any $X,Y \in \Ind(N(\bfT))$. Since $F$ preserves filtered $\infty$-colimits, it suffices to show this for $X \in N(\bfT)$. Write $Y = \colim_i Y_i$ as a filtered $\infty$-colimit of a diagram $Y \colon I \to N(\bfT)$ (which we also denote by $Y$). By Proposition 5.3.1.18 of \cite{Lurie2009HTT} and Lemma E.1.6.4 of \cite{lurie2016spectral}, we may assume without loss of generality that $I$ is the nerve of a directed poset, which we also denote by $I$, with the property that $I_{< i}$ is finite for any $i \in I$. By \autoref{lem:Rectification}, we may replace $Y$ by a strict diagram $Z \colon I \to \bfT$. Since a diagram as described in \autoref{lem:Rectification} is cofibrant in the projective model structure on $\Ind(\bfC)^I$, we see that the ind-object $Z = \{Z_i\}_{i \in I}$ is the homotopy colimit of the diagram $i \mapsto Z_i$. By Theorem 4.2.4.1 of \cite{Lurie2009HTT}, the object $Z$ is an $\infty$-colimit of the diagram $Y \colon N(I) \to \Ind(\bfC)_\infty$, hence $Z$ is equivalent to $F(Y)$ (note that $F$ preserves filtered colimits). In particular, we obtain a commutative diagram
\[\begin{tikzcd}
\colim_i \map_{\Ind(N(\bfT))}(X,Z_i) \ar[r] \ar[d] & \colim_i \map_{\Ind(\bfC)_\infty}(FX,FZ_i) \ar[d]\\
\map_{\Ind(N(\bfT))}(X,Y) \ar[r] & \map_{\Ind(\bfC)_\infty}(FX,FY).
\end{tikzcd}\]
Here the left-hand vertical map is an equivalence since objects of $\bfT$ are compact (in the $\infty$-categorical sense), while the right-hand vertical map is an equivalence since it is equivalent to $\colim_i \Map_{\Ind(\bfC)}(X,Z_i) \to \Map_{\Ind(\bfC)}(X,Z)$, which is an isomorphism since $X$ is compact in $\Ind(\bfC)$. The top horizontal map is an equivalence since $F$ is by construction fully faithful when restricted to $N(\bfT) \subset \Ind(N(\bfT))$. We conclude that the bottom map is an equivalence and hence that $F$ is fully faithful.

To see that $F$ is essentially surjective, let $X$ be a fibrant-cofibrant object in $\Ind(\bfC)$. By \autoref{lemma:CofibrantObjectIsDirectedColimitTestObjects}, $X$ is a directed colimit $\colim_i t_i$ of objects in $\bfT$. By \autoref{lemma:WeakEquivalencesStableUnderFilCol}, $X$ is also a homotopy colimit of this diagram, hence $X$ is an $\infty$-limit of the diagram $\{t_i\}_i$ in the underlying $\infty$-category $\Ind(\bfC)_\infty$. View $\{t_i\}_i$ as a diagram in $N(\bfT)$ and let $Y$ denote the $\infty$-colimit of this diagram in $\Ind(N(\bfT))$. Since $F$ preserves filtered $\infty$-colimits, it follows that $F(Y) \simeq X$ and hence that $F$ is essentially surjective.
\end{proof}

We automatically obtain the following dual result. Note that item \eqref{thm:ApproximateTheoremItem4} of \autoref{thm:ApproximateTheorem} stated in the introduction is a direct consequence of this theorem.

\begin{theorem}\label{thm:UnderlyingInftyCategoryPro}
Let $(\bfC,\bfT)$ be a fibration test category and suppose that $\bfT$ is closed under pullbacks along fibrations (see \autoref{definition:ClosedUnderPullbackAlongCofib}). Then the canonical functor
\[\Pro(N(\bfT)) \to \Pro(\bfC)_\infty\]
is an equivalence of quasi-categories.
\end{theorem}

The main theorems of this appendix can be used to determine the underlying $\infty$-categories of many of the examples that were mentioned throughout this paper. Moreover, it shows that the homotopy theory of $\Pro(\bfC)$ is often fully determined by the full simplicial subcategory $\bfT$ of $\bfC$. By way of illustration, we will to single out one specific example. Namely, we will relate the ``profinite'' Joyal-Kan model structure (see \autoref{example:PStratifiedFibrationTestCategory}) to the profinite stratified spaces defined in \cite[\S 2.5]{BarwickGlasmanHaine2018ExodromyV7}. Note that one can use similar arguments to determine the underlying $\infty$-categories of Quick's and Morel's model structures on $\s\wh\Set$ (cf.~\cite[\S 7]{BarneaHarpazHorel2017}) and of the profinite Joyal model structure (see \autoref{remark:UnderlyingInftyCategoryProfiniteInftyCats}).

\begin{example}\label{example:UnderlyingInftyCatProfiniteStratified}
Let $P$ be a finite poset and let $\LsSet_{/P}$ be the fibration test category defined in \autoref{example:PStratifiedFibrationTestCategory}. The full subcategory of test objects $\bfT$ in this fibration test category consists of the fibrant objects of the Joyal-Kan model structure on $\s\Set_{/P}$ whose total space is a lean simplicial set. We will call these lean $P$-stratified Kan complexes. They can be described explicitly as those inner fibrations $f \colon X \fibarrow P$ for which $X$ is lean and the fiber above any point is a Kan complex. It is proved in \autoref{lem:LeanStratifiedKanComplexes} below that the homotopy-coherent nerve $N(\bfT)$ of the category of lean $P$-stratified Kan complexes is equivalent to $\mathbf{Str}_{\pi,P}$, the $\infty$-category of $\pi$-finite $P$-stratified spaces defined in Definition 2.4.3 of \cite{BarwickGlasmanHaine2018ExodromyV7}. By \autoref{thm:UnderlyingInftyCategoryPro}, it now follows that the underlying $\infty$-category of the profinite Joyal-Kan model structure on $\s\wh\Set_{/P}$ is equivalent to $\Pro(\mathbf{Str}_{\pi,P})$, which is equivalent to the $\infty$-category of profinite $P$-stratified spaces defined in \cite[\S 2.5]{BarwickGlasmanHaine2018ExodromyV7}.
\end{example}

We conclude this appendix by proving the lemma used in the above example.

\begin{lemma}\label{lem:LeanStratifiedKanComplexes}
Let $P$ be (the nerve of) a finite poset and let $\bfT$ be the full simplicial subcategory of $\s\Set_{/P}$ spanned by the lean $P$-stratified Kan complexes. Then the homotopy-coherent nerve $N(\bfT)$ is equivalent to the $\infty$-category of $\pi$-finite $P$-stratified spaces as defined in Definition 2.4.3 of \cite{BarwickGlasmanHaine2018ExodromyV7}.
\end{lemma}

\begin{proof}
By slightly rephrasing the definition of ``$\pi$-finite''  given in \cite{BarwickGlasmanHaine2018ExodromyV7}, this comes down to proving that if $X \fibarrow P$ is a lean $P$-stratified Kan complex, then
\begin{enumerate}[(i)]
    \item\label{item1pifinite} for any $p \in P$, the set $\pi_0(f^{-1}(p))$ is finite,
    \item\label{item2pifinite} there exists an $n \in \bbN$ such that, for all $x,y \in X$, the homotopy groups of $\map_X(x,y)$ vanish above degree $n$, and
    \item\label{item3pifinite} for all $x,y \in X$, the Kan complex $\map_X(x,y)$ has finite homotopy groups
\end{enumerate}
and conversely that any $P$-stratified Kan complex $X \fibarrow P$ satisfying these properties is equivalent to a lean $P$-stratified Kan complex. If $X$ is a lean, then items \eqref{item1pifinite} and \eqref{item3pifinite} follows since $X$ is degreewise finite, while \eqref{item2pifinite} follows since $X$ is coskeletal. For the converse, let a $P$-stratified Kan complex $X \fibarrow P$ satisfy these items. If we replace $X \fibarrow P$ by a minimal inner fibration $\wt X \to P$ (cf.~\cite[\S 2.3.3]{Lurie2009HTT}), then it is still a $P$-stratified Kan complex satisfying items \eqref{item1pifinite}-\eqref{item3pifinite}, so it suffices to show that $\wt X$ is lean. Since pullbacks of minimal fibrations are again minimal, it follows from \eqref{item1pifinite} that $f^{-1}(p) \subset \wt X$ has finitely many $0$-simplices for any $p \in P$, and hence that $\wt X$ has finitely many $0$-simplices. Since $P$ is (the nerve of) a poset, two maps $\Delta^n \to \wt X$ are homotopic relative to the boundary if and only if they are so over $P$. This implies that $\wt X$ is itself a minimal quasi-category, and hence degreewise finite by \eqref{item3pifinite} and \autoref{lem:MinimalPiFiniteQuasiCategoryIsDegreewiseFinite} below. It is proved in Proposition 2.3.4.18 of \cite{Lurie2009HTT} that if $\wt X$ is a minimal quasi-category satisfying \eqref{item2pifinite}, then it is coskeletal, so we conclude that $\wt X$ is lean.
\end{proof}

\begin{lemma}\label{lem:MinimalPiFiniteQuasiCategoryIsDegreewiseFinite}
Let $X$ be a minimal quasi-category with finitely many $0$-simplices and with the property that for any $x,y \in X_0$, the homotopy groups of $\map_X(x,y)$ are finite. Then $X$ is degreewise finite.
\end{lemma}

\begin{proof}
Since $X$ has finitely many $0$-simplices, it suffices to show that for any $n \geq 1$ and any map $D \colon \partial \Delta^n \to X$, there exist finitely many $n$-simplices filling $D$. For $n=1$ this is clear: by minimality, the number of $1$-simplices from $x$ to $y$ in $X$ agrees with $\pi_0 \map(x,y)$, which is finite by assumption. Now assume $n > 1$ and let $D$ be given. Write $E$ for the restriction of $D$ to the face opposite to the $n$-th vertex, and write $\partial E$ for the restriction of $E$ to $\partial \Delta^{n-1}$. This restriction induces a left fibration $X_{E/} \to X_{\partial E /}$, where these slice categories are defined as in \cite[\S 3]{Joyal2002QuasiKan}. Let $z$ be the $0$-simplex of $X$ obtained by restricting $D$ to the top vertex, and denote the fibers of $X_{E/}$ and $X_{\partial E /}$ above $z$ by $\map(E,z)$ and $\map(\partial E, z)$, respectively. Since these are fibers of left fibrations over $X$, we see that these are Kan complexes. Note that the restriction of $D$ to $\Lambda^n_n$ defines a $0$-simplex in $\map(\partial E, z)$. Now define $\mathrm{Fill}(D)$ as the pullback
\[\begin{tikzcd}
\mathrm{Fill}(D) \arrow[dr,phantom, very near start, "\lrcorner"] \ar[r] \ar[d] & \map(E,z) \ar[d, two heads] \\
\{D|_{\Lambda^n_n}\} \ar[r] & \map(\partial E, z).
\end{tikzcd}\]
It is clear that the $0$-simplices of $\mathrm{Fill}(D)$ correspond to $n$-simplices in $X$ that fill $D$. A $1$-simplex in $\mathrm{Fill}(D)$ between two such $n$-simplices $f,g$ in $X$ is exactly an $(n+1)$-simplex $h \colon \Delta^{n+1} \to X$ such that $d_n h = f$, $d_{n+1} h = g$ and $d_i h = d_i s_m f$ for any $i < n$. Given such an $(n+1)$-simplex $h$, the sequence $(s_0 f, s_1 f, \ldots, s_{n-1} f, h)$ defines a homotopy $\Delta^n \times \Delta^1 \to X$ between $f$ and $g$ relative to $\partial \Delta^n$.\footnote{In fact, the converse is also true: if there is a homotopy between $f$ and $f'$ relative to $\partial \Delta^n$, then there exists an $(n+1)$-simplex $h$ in $X$ with the given property. A proof of this statement can be obtained by slightly modifying the proof of Theorem I.8.2 in \cite{Lamotke1968Semisimpliziale} in such a way that one only needs to fill inner horns.} In particular, by minimality of $X$, the existence of such an $(n+1)$-simplex $h$ implies that $f = g$, and hence the number of elements in $\pi_0(\mathrm{Fill}(D))$ equals the number of fillers of $D \colon \partial \Delta^n \to X$.

Since $X_{E/} \to X_{\partial E /}$ is a left fibration and $\map(\partial E, z)$ is a Kan complex, the restriction $\map(E,z) \to \map(\partial E, z)$ is a Kan fibration and hence $\mathrm{Fill}(D)$ is the homotopy fiber of $\map(E,z) \fibarrow \map(\partial E, z)$. In particular, if $\map(E,z)$ and $\map(\partial E, z)$ have finite homotopy groups, then $\mathrm{Fill}(D)$ does as well, and hence $D$ has finitely many fillers. If we let $y$ denote the top vertex of the $(n-1)$-simplex $E$, then $\map(E,z) \simeq \map(y,z)$, which has finite homotopy groups by assumption. To see that $\map(\partial E, z)$ has finite homotopy groups, note that
\[\map(\partial E, z) = \lim_{x \in \mathrm{nd}(\partial \Delta^n)^{op}} \map(E|_x, z)\]
where $\mathrm{nd}(\partial \Delta^n)$ denotes the poset of non-degenerate simplices of $\partial \Delta^{n-1}$. This follows from the fact that the join of simplicial sets $\star$ preserves connected colimits. We see that for any $x \in \mathrm{nd}(\partial \Delta^n)$, the Kan complex $\map(E|_x,z)$ is equivalent to $\map(y,z)$, where $y$ denotes the top vertex of $E|_x$. In particular, it has finite homotopy groups. Note that the diagram $x \mapsto \map(E|_x,z)$ is injectively fibrant since the diagram $\{x\}_{x \in \mathrm{nd}(\partial \Delta^n)}$ is cofibrant in the projective model structure on $\s\Set^{\mathrm{nd}(\partial \Delta^n)}$. In particular, $\map(\partial E, z)$ is a finite homotopy limit of spaces with finite homotopy groups, so it has finite homotopy groups as well. We conclude that there are finitely many $n$-simplices filling $D \colon \partial \Delta^n \to X$.
\end{proof}

\begin{remark}\label{remark:UnderlyingInftyCategoryProfiniteInftyCats}
It follows as in the proofs of \Cref{lem:LeanStratifiedKanComplexes,lem:MinimalPiFiniteQuasiCategoryIsDegreewiseFinite} that a quasi-category is equivalent to a lean quasi-category if and only it has finitely many objects up to equivalence and all its mapping spaces have finite homotopy groups that vanish above a certain dimension; let us call such quasi-categories \emph{$\pi$-finite}. Applying \autoref{thm:UnderlyingInftyCategoryPro} to the fibration test category $\LsSet_J$ of \autoref{example:JoyalFibrationTestCategory} shows that the underlying $\infty$-category of the profinite Joyal model structure $\s\wh\Set_J$ (and hence also of $\biswhSet_{CSS}$) is equivalent to $\Pro(\Cat_{\infty, \pi})$, where $\Cat_{\infty, \pi}$ denotes the $\infty$-category of $\pi$-finite $\infty$-categories.
\end{remark}

\printbibliography

\noindent{\sc Matematiska institutionen, Stockholms universitet, 106 91 Stockholm, Sweden.}\\
\noindent{\emph{E-mail:} \href{mailto:blom@math.su.se}{\nolinkurl{blom@math.su.se}}}\vspace{2ex}

\noindent{\sc Mathematisch Instituut, Universiteit Utrecht, Postbus 80010, 3508 TA Utrecht, The Netherlands.}\\
\noindent{\emph{E-mail:} \href{mailto:i.moerdijk@uu.nl}{\nolinkurl{i.moerdijk@uu.nl}}}

\end{document}